\documentclass[12pt,oneside,reqno]{amsart}
\usepackage[margin=1in]{geometry}
\usepackage[utf8]{inputenc}
\usepackage{amsthm, amssymb}
\usepackage[colorlinks=true, pdfstartview=FitV, linkcolor=blue, citecolor=blue, urlcolor=blue]{hyperref}
\usepackage[colorinlistoftodos]{todonotes}
\usepackage{algorithm}
\usepackage{algpseudocode}
\usepackage{pgf}
\usepackage{tikz}
\usepackage{tikz-cd}
\usetikzlibrary{positioning,shapes,shadows,arrows}
\usepackage{bbm,bm}
\usepackage{enumitem}
\usepackage{mathabx}
\usepackage{comment}
\usepackage{ytableau}
\usepackage{thmtools}
\usepackage{thm-restate}

\newtheorem{prop}{Proposition}[section]
\newtheorem{thm}[prop]{Theorem}
\newtheorem{lem}[prop]{Lemma}
\newtheorem{lemma}[prop]{Lemma}
\newtheorem{cor}[prop]{Corollary}

\theoremstyle{remark}
\newtheorem{exa}[prop]{Example}
\newtheorem{rem}[prop]{Remark}
\newtheorem{defn}[prop]{Definition}

\newcommand{\Z}{\mathbb{Z}}

\newcommand{\fG}{\mathfrak{G}}

\newcommand{\IR}{\mathsf{IR}} 

\newcommand{\wt}{\mathsf{wt}}

\newcommand{\rajcode}{\mathsf{rajcode}}

\newcommand{\reg}{\mathsf{reg}}

\newcommand{\Rothe}{\mathsf{Rothe}}
\newcommand{\PD}{\mathsf{PD}}
\newcommand{\dark}{\mathsf{dark}}
\newcommand{\invcode}{\mathsf{invcode}}
\newcommand{\movecode}{\mathsf{movecode}}

\newcommand{\jtile}{
\begin{tikzpicture}[x=1em,y=1em,thick,rounded corners,color = red]
\draw[step=1,gray,thin] (0,0) grid (1,1);
\draw[color=black, thick, sharp corners] (0,0) rectangle (1,1);
\draw(0.5,1.0)--(0.5,0.5)--(0.0,0.5);
\end{tikzpicture}}

\newcommand{\bumptile}{
\begin{tikzpicture}[x=1em,y=1em,thick,rounded corners,color = red]
\draw[step=1,gray,thin] (0,0) grid (1,1);
\draw[color=black, thick, sharp corners] (0,0) rectangle (1,1);
\draw(1.0,0.5)--(0.5,0.5)--(0.5,0.0);
\draw(0.5,1.0)--(0.5,0.5)--(0.0,0.5);
\end{tikzpicture}}

\newcommand{\ptile}{
\begin{tikzpicture}[x=1em,y=1em,thick,color = red]
\draw[step=1,gray,thin] (0,0) grid (1,1);
\draw[color=black, thick, sharp corners] (0,0) rectangle (1,1);
\draw(0.5,1.0)--(0.5,0.0);
\draw(1.0,0.5)--(0.0,0.5);
\end{tikzpicture}}

\definecolor{darkblue}{rgb}{0.0,0,0.7} 
\definecolor{darkred}{rgb}{0.7,0,0} 
\definecolor{darkgreen}{rgb}{0, .6, 0} 
\newcommand{\definition}[1]{{\color{darkred}\emph{#1}}} 

\title{Constructing maximal pipedreams of double Grothendieck polynomials}

\author[C.~Chou]{Chen-An Chou}
\address[C. Chou]{Department of Mathematics, UC San Diego, La Jolla, CA 92093, U.S.A.}
\email{c1chou@ucsd.edu}

\author[T.~Yu]{Tianyi Yu}
\address[T. Yu]{Department of Mathematics, UC San Diego, La Jolla, CA 92093, U.S.A.}
\email{tiy059@ucsd.edu}

\begin{document}
\maketitle
\begin{abstract}
Pechenik, Speyer and Weigandt 
defined a statistic $\rajcode(\cdot)$
on permutations which characterizes
the leading monomial in top degree components
of double Grothendieck polynomials.
Their proof is combinatorial: 
They showed there
exists a unique pipedream
of a permutation $w$ with row weight $\rajcode(w)$
and column weight $\rajcode(w^{-1})$.
They proposed the problem of finding a ``direct recipe''  for this pipedream.
We solve this problem 
by providing an algorithm
that constructs this pipedream via ladder moves.
\end{abstract}
\section{Introduction}
\label{S: Intro}

The matrix Schubert variety $X_w$ is a determinantal variety
that has been studied extensively
(see for instance~\cite{Ful, KM, KMY, WY}).
Castelnuovo–Mumford regularity measures the algebraic complexity of varieties.
Since matrix Schubert varieties are Cohen–Macaulay~\cite{Ful, KM, Ram}, the Castelnuovo–Mumford regularity of $X_w$ is the difference between the top and bottom degree of its K-polynomial.
By the work of Knutson and Miller~\cite{KM}, the K-polynomial of $X_w$ is the Grothendieck polynomial $\fG_w(\textbf{x})$.
This family of polynomials, introduced by Lascoux and Sch\"utzenberger \cite{LS:Groth}, represents K-classes of structure sheaves of Schubert varieties in flag varieties. 
Their lowest degree components are the Schubert polynomials whose degrees are known.

Consequently, determining the Castelnuovo–Mumford regularity of $X_w$ reduces to computing the degree of $\fG_w(\textbf{x})$.
With this motivation, there has been a recent surge in the study of top degree components of $\fG_w(\textbf{x})$~\cite{DMS, Haf, PSW, PY, RRRSW, RRW}.
Pechenik, Speyer, and Weigandt~\cite{PSW} defined a statistic $\rajcode(\cdot)$ on $S_n$
using increasing subsequences of permutations.
They showed $x^{\rajcode(w)}$
is the leading monomial in the top degree
components of $\fG_w(\textbf{x})$
with respect to the lexicographical order 
where $x_n > \cdots > x_1$.
Pan and Yu~\cite{PY} found a diagrammatic formula to compute $\rajcode(w)$ 
(see Definition~\ref{D: rajcode}).

For $w \in X_n$, the double Grothendieck polynomial $\fG_w(\textbf{x},\textbf{y})$ 
involves two sets of variables: $x_1, \cdots, x_n$ and $y_1, \cdots, y_n$.
It represents Schubert classes in the torus-equivariant K-theory of flag varieties.
After setting $y_1 = y_2 = \cdots = 0$, the double Grothendieck polynomial
$\fG_w(\textbf{x},\textbf{y})$ specializes to the usual Grothendieck polynomial $\fG_w(\textbf{x})$.
The $\rajcode(\cdot)$ statistic 
also captures the leading monomial
in top degree components of $\fG_w(\textbf{x},\textbf{y})$.

\begin{thm}[{\cite[Theorem 1.4]{PSW}}]
\label{T: PSW}
The leading monomial of top degree components of $\fG_w(\textbf{x},\textbf{y})$ is $x^{\rajcode(w)} y^{\rajcode(w^{-1})}$ with coefficient $1$
for any term order with 
$x_n > \cdots > x_1$
and $y_n > \cdots > y_1$.
\end{thm}

A combinatorial formula of $\fG_w(\textbf{x},\textbf{y})$ is given by pipedreams~\cite{BB, BJS, FK, KM}: certain tilings of a staircase grid 
using $\ptile$, $\bumptile$ and $\jtile$ (see Definition~\ref{D: PD}).
The row (resp. column) weight of a pipedream is a weak composition where the $i^\textsuperscript{th}$
entry is the number of $\ptile$ in row (resp. column) $i$ of the pipedream.
Let $\PD(w)$ be the set of pipedreams 
associated with the permutation $w$.
Pechenik, Speyer, and Weigandt established Theorem~\ref{T: PSW} by showing there exists a unique pipedream in $\PD(w)$ with row weight $\rajcode(w)$ and column weight $\rajcode(w^{-1})$,
which they call the maximal pipedream of $w$. 
In Remark 7.2, they said:
\begin{quote}
``We find it frustrating that we do not have a direct recipe for the maximal
pipe dream in terms of 
$w$.''    
\end{quote}
The main goal of this paper is to relieve their frustration:
We give an explicit algorithm to construct the maximal pipedream 
$\widehat{P}(w) \in \PD(w)$.

\begin{thm}
\label{T: Row col weight}
For $w \in S_n$,
the pipedream $\widehat{P}(w)$ we construct has row weight $\rajcode(w)$ and column weight $\rajcode(w^{-1})$.
\end{thm}

Our algorithm involves a local move
known as the ladder move~\cite{BB}.
When row $r$ column $c$ of a pipedream $P$
is $\ptile$, we write $(r,c) \in P$. 
We may apply a ladder move on a $\ptile$ 
in row $r$ column $c$ of a pipedream $P$ if
all the following are satisfied: 
\begin{itemize}
\item $(r, c+1) \notin P$.
\item There exists $r' < r$ such that 
$(r',c) \notin P$ and $(r',c + 1) \notin P$.
In addition, $(i,c), (i,c+1) \in P$
for any $r' < i < r$.
\end{itemize}
Now we perform the ladder move at the 
$\ptile$ in row $r$ column $c$ of $P$.
First turn the $\bumptile$ at row $r'$ 
column $c+1$ into a $\ptile$.
Then we may or may not turn 
the $\ptile$ at row $r$ 
column $c$ into $\bumptile$.
If we do that,
the move is called a \definition{regular ladder move}.
Otherwise, the move is called a 
\definition{K-ladder move}.
Locally, the moves look like the following:

$$
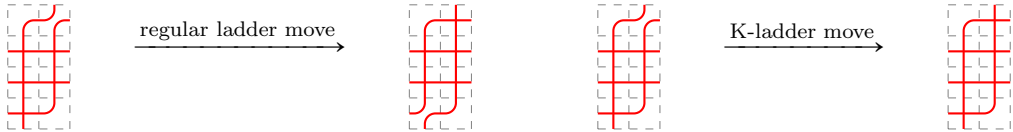
\begin{tikzpicture}[x=1em,y=1em,thick, rounded corners,color = red]
\draw[step=1,gray,ultra thin,dashed] (0,0) grid (2,4);
\draw[color=red] (0,0.5)--(1.5, 0.5)--(1.5, 3.5)--(2, 3.5);
\draw[color=red] (0.5,0)--(0.5, 3.5)--(1.5, 3.5)--(1.5, 4);
\draw[color=red] (0,1.5)--(2, 1.5);
\draw[color=red] (0,2.5)--(2, 2.5);
\end{tikzpicture}
\quad\quad \raisebox{1cm}{$\xrightarrow{\textrm{regular ladder move}}$}\quad\quad
\begin{tikzpicture}[x=1em,y=1em,thick, rounded corners,color = red]
\draw[step=1,gray,ultra thin,dashed] (0,0) grid (2,4);
\draw[color=red] (0,0.5)--(0.5, 0.5)--(0.5, 3.5)--(2, 3.5);
\draw[color=red] (0.5,0)--(0.5, 0.5)--(1.5, 0.5)--(1.5, 4);
\draw[color=red] (0,1.5)--(2, 1.5);
\draw[color=red] (0,2.5)--(2, 2.5);
\end{tikzpicture}
\quad\quad\quad\quad
\begin{tikzpicture}[x=1em,y=1em,thick, rounded corners,color = red]
\draw[step=1,gray,ultra thin,dashed] (0,0) grid (2,4);
\draw[color=red] (0,0.5)--(1.5, 0.5)--(1.5, 3.5)--(2, 3.5);
\draw[color=red] (0.5,0)--(0.5, 3.5)--(1.5, 3.5)--(1.5, 4);
\draw[color=red] (0,1.5)--(2, 1.5);
\draw[color=red] (0,2.5)--(2, 2.5);
\end{tikzpicture}
\quad\quad\raisebox{1cm}{$\xrightarrow{\textrm{K-ladder move}}$}\quad\quad
\begin{tikzpicture}[x=1em,y=1em,thick, rounded corners,color = red]
\draw[step=1,gray,ultra thin,dashed] (0,0) grid (2,4);
\draw[color=red] (0,0.5)--(1.5, 0.5)--(1.5, 4);
\draw[color=red] (0.5,0)--(0.5, 3.5)--(2, 3.5);
\draw[color=red] (0,1.5)--(2, 1.5);
\draw[color=red] (0,2.5)--(2, 2.5);
\end{tikzpicture}
$$

For $w \in S_n$,
the statistic $\invcode(w)$
is a sequence of $n$ numbers where
the $i\textsuperscript{th}$ number
is the number of $j > i$ such that
$w(j) < w(i)$.
It is well-known that $\PD(w)$
contains the pipedream
with row weight $\invcode(w)$
and all $\ptile$ are left-justified.
All other pipedreams in $\PD(w)$
can be obtained by performing
ladder moves from this one. 
We start from this pipedream
and perform an iterative algorithm.
Each iteration places a bar right above row $i$ for $i = n-2, n-3, \cdots, 1$.
During each iteration, 
we only look under the bar and imagine row $i$ is the topmost row. 
Scan through the columns from right to left. 
Within each column, scan through the $\ptile$ from top to bottom. 
Whenever we see a $\ptile$ at which we can perform a ladder move, we perform a regular ladder move. 
After going through a column, if we have performed ladder moves on this column, we turn the last ladder move into a K-ladder move.
We denote the final pipedream by $\widehat{P}(w)$.

\begin{exa}
\label{E: Algorithm}
Take $w \in S_5$ with one-line notation $14523$.
We start from the following pipedream:
$$
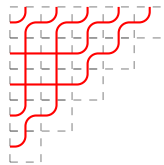
\begin{tikzpicture}[x=1em,y=1em,thick,rounded corners, color = blue]
\draw[step=1,gray,ultra thin,dashed] (0,0) grid (1,1);
\draw[step=1,gray,ultra thin,dashed] (0,1) grid (2,2);
\draw[step=1,gray,ultra thin,dashed] (0,2) grid (3,3);
\draw[step=1,gray,ultra thin,dashed] (0,3) grid (4,4);
\draw[step=1,gray,ultra thin,dashed] (0,4) grid (5,5);
\draw[color=red] (0,0.5)--(0.5, 0.5)--(0.5, 1.5)--(1.5, 1.5)--(1.5, 4.5)--(2.5, 4.5)--(2.5, 5);
\draw[color=red] (0,1.5)--(0.5, 1.5)--(0.5, 4.5)--(1.5, 4.5)--(1.5, 5);
\draw[color=red] (0,2.5)--(2.5, 2.5)--(2.5, 3.5)--(3.5, 3.5)--(3.5, 4.5)--(4.5, 4.5)--(4.5, 5);
\draw[color=red] (0,3.5)--(2.5, 3.5)--(2.5, 4.5)--(3.5, 4.5)--(3.5, 5);
\draw[color=red] (0,4.5)--(0.5, 4.5)--(0.5, 5);
\end{tikzpicture}
$$
When $i = 3$ and $2$, we do not make any moves. 
When $i = 1$, we perform:
$$
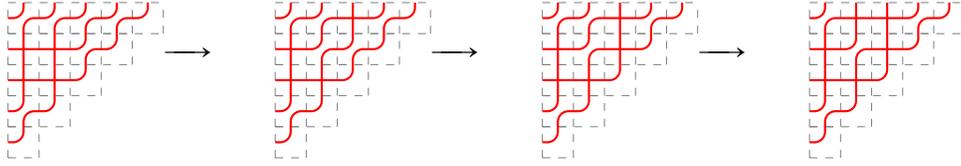
\begin{tikzpicture}[x=1em,y=1em,thick,rounded corners, color = blue]
\draw[step=1,gray,ultra thin,dashed] (0,0) grid (1,1);
\draw[step=1,gray,ultra thin,dashed] (0,1) grid (2,2);
\draw[step=1,gray,ultra thin,dashed] (0,2) grid (3,3);
\draw[step=1,gray,ultra thin,dashed] (0,3) grid (4,4);
\draw[step=1,gray,ultra thin,dashed] (0,4) grid (5,5);
\draw[color=red] (0,0.5)--(0.5, 0.5)--(0.5, 1.5)--(1.5, 1.5)--(1.5, 4.5)--(2.5, 4.5)--(2.5, 5);
\draw[color=red] (0,1.5)--(0.5, 1.5)--(0.5, 4.5)--(1.5, 4.5)--(1.5, 5);
\draw[color=red] (0,2.5)--(2.5, 2.5)--(2.5, 3.5)--(3.5, 3.5)--(3.5, 4.5)--(4.5, 4.5)--(4.5, 5);
\draw[color=red] (0,3.5)--(2.5, 3.5)--(2.5, 4.5)--(3.5, 4.5)--(3.5, 5);
\draw[color=red] (0,4.5)--(0.5, 4.5)--(0.5, 5);
\end{tikzpicture}
 \raisebox{1.3cm}{$\xrightarrow{\quad}$}\quad\quad
\begin{tikzpicture}[x=1em,y=1em,thick,rounded corners, color = blue]
\draw[step=1,gray,ultra thin,dashed] (0,0) grid (1,1);
\draw[step=1,gray,ultra thin,dashed] (0,1) grid (2,2);
\draw[step=1,gray,ultra thin,dashed] (0,2) grid (3,3);
\draw[step=1,gray,ultra thin,dashed] (0,3) grid (4,4);
\draw[step=1,gray,ultra thin,dashed] (0,4) grid (5,5);
\draw[color=red] (0,0.5)--(0.5, 0.5)--(0.5, 1.5)--(1.5, 1.5)--(1.5, 3.5)--(2.5, 3.5)--(2.5, 5);
\draw[color=red] (0,1.5)--(0.5, 1.5)--(0.5, 4.5)--(1.5, 4.5)--(1.5, 5);
\draw[color=red] (0,2.5)--(2.5, 2.5)--(2.5, 3.5)--(3.5, 3.5)--(3.5, 4.5)--(4.5, 4.5)--(4.5, 5);
\draw[color=red] (0,3.5)--(1.5, 3.5)--(1.5, 4.5)--(3.5, 4.5)--(3.5, 5);
\draw[color=red] (0,4.5)--(0.5, 4.5)--(0.5, 5);
\end{tikzpicture}
 \raisebox{1.3cm}{$\xrightarrow{\quad}$}\quad\quad
\begin{tikzpicture}[x=1em,y=1em,thick,rounded corners, color = blue]
\draw[step=1,gray,ultra thin,dashed] (0,0) grid (1,1);
\draw[step=1,gray,ultra thin,dashed] (0,1) grid (2,2);
\draw[step=1,gray,ultra thin,dashed] (0,2) grid (3,3);
\draw[step=1,gray,ultra thin,dashed] (0,3) grid (4,4);
\draw[step=1,gray,ultra thin,dashed] (0,4) grid (5,5);
\draw[color=red] (0,0.5)--(0.5, 0.5)--(0.5, 1.5)--(1.5, 1.5)--(1.5, 3.5)--(3.5, 3.5)--(3.5, 4.5)--(4.5, 4.5)--(4.5, 5);
\draw[color=red] (0,1.5)--(0.5, 1.5)--(0.5, 4.5)--(1.5, 4.5)--(1.5, 5);
\draw[color=red] (0,2.5)--(2.5, 2.5)--(2.5, 5);
\draw[color=red] (0,3.5)--(1.5, 3.5)--(1.5, 4.5)--(3.5, 4.5)--(3.5, 5);
\draw[color=red] (0,4.5)--(0.5, 4.5)--(0.5, 5);
\end{tikzpicture}
\raisebox{1.3cm}{$\xrightarrow{\quad}$}\quad\quad
\begin{tikzpicture}[x=1em,y=1em,thick,rounded corners, color = blue]
\draw[step=1,gray,ultra thin,dashed] (0,0) grid (1,1);
\draw[step=1,gray,ultra thin,dashed] (0,1) grid (2,2);
\draw[step=1,gray,ultra thin,dashed] (0,2) grid (3,3);
\draw[step=1,gray,ultra thin,dashed] (0,3) grid (4,4);
\draw[step=1,gray,ultra thin,dashed] (0,4) grid (5,5);
\draw[color=red] (0,0.5)--(0.5, 0.5)--(0.5, 1.5)--(1.5, 1.5)--(1.5, 3.5)--(3.5, 3.5)--(3.5, 4.5)--(4.5, 4.5)--(4.5, 5);
\draw[color=red] (0,1.5)--(0.5, 1.5)--(0.5, 4.5)--(3.5, 4.5)--(3.5, 5);
\draw[color=red] (0,2.5)--(2.5, 2.5)--(2.5, 5);
\draw[color=red] (0,3.5)--(1.5, 3.5)--(1.5, 5);
\draw[color=red] (0,4.5)--(0.5, 4.5)--(0.5, 5);
\end{tikzpicture}
$$
\end{exa}

Dreyer, M{\'e}sz{\'a}ros, and St. Dizier~\cite{DMS}
found the leading monomial in each homogeneous
component of $\fG_w$.
Let $\reg(w)$ be the difference
between the sum of entries in $\rajcode(w)$
and the sum of entries in $\invcode(w)$.
Define the map $\IR(\cdot)$
that sends $w$ to a sequence of monomials
$m_0, m_1, \cdots, m_{\reg(w)}$.
First, $m_0 := x^{\invcode(w)}$.
For $i > 0$, $m_i := m_{i-1} x_p$
where $p$ is the largest such that 
$m_{i-1} x_p$ divides $x^{\rajcode(w)}$.
For each $m_i$, 
Dreyer, M{\'e}sz{\'a}ros, and St. Dizier~\cite{DMS} explicitly constructed a climbing chain,
another combinatorial model of $\fG_w$ introduced in~\cite{LRS},
showing $m_i$ is the leading monomial
in its degree of $\fG_w$.
In our algorithm, we start from a pipedream
with row weight $\invcode(w)$.
During the algorithm, 
we obtain the pipedreams corresponding to 
$m_1,
\cdots, m_{\reg(w)}$.  

\begin{thm}
\label{T: Intermidiate row}
Let $w \in S_n$.
Perform our algorithm to compute $\widehat{P}(w)$.
The algorithm makes $\reg(w)$
K-ladder moves. 
Right after the $i^\textsuperscript{th}$
K-ladder move,
we record the row weight of the pipedream as $a_i(w)$.
Then $x^{a_i(w)} = m_i$
where $\IR(w) = (m_0, m_1, \cdots, m_{\reg(w)})$.
\end{thm}

The rest of the paper is structured as follows.
In~\S\ref{S: Background},
we cover necessary background 
regarding pipedreams and $\rajcode(w)$.
In~\S\ref{S: Recursion},
we introduce recursive 
formulas to compute $\rajcode(w)$, 
$\rajcode(w^{-1})$ and $\IR(w)$.
In~\S\ref{S: Proofs},
we prove our main results
using Proposition~\ref{P: Movecode}
and Corollary~\ref{C: Row +1}, 
whose proofs are in~\S\ref{S: Proofs of Prop and Cor}.

\section{Background}
\label{S: Background}

\subsection{Pipedreams and Grothendieck polynomials}
\begin{defn}
\label{D: PD}
\definition{Pipedreams} of size $n$ 
are tilings with 
$n + 1 - i$ left justified tiles in row $i$.
The rightmost tile in each row is $\jtile$
and all other tiles can be $\bumptile$
or $\ptile$.
For a pipedream of size $n$,
it is associated with a permutation $w \in S_n$.
We label the pipes $1,2, \cdots, n$ along the
top edge and follow the pipes. 
Whenever two pipes cross more than once, 
we treat all but the first crossing
as $\bumptile$.
Let $\PD(w)$ be the set of the pipedreams
associated with $w \in S_n$.
\end{defn}

\begin{exa}
Pipedreams in Example~\ref{E: Algorithm}
are all in $\PD(w)$ where $w$ has one-line
notation $14523$.
\end{exa}

Let $P$ be a pipedream. 
We write $(i,j) \in P$
if row $i$ column $j$ of $P$
is $\ptile$.
Following~\cite{KM} and~\cite{FK},
\definition{double Grothendieck polynomial} $\fG_w(\textbf{x}, \textbf{y})$
and \definition{Grothendieck polynomial} $\fG_w(\textbf{x})$
can be defined as 
$$
\fG_w(\textbf{x}, \textbf{y}) := \sum_{P \in \PD(w)}\prod_{(i,j) \in P} (x_i + y_j - x_i y_j),
\quad\quad
\fG_w(\textbf{x}) := \sum_{P \in \PD(w)}\prod_{(i,j) \in P} x_i.
$$

In the rest of the paper, 
we identify a pipedream with a \definition{diagram}, which is a finite subset of
$\Z_{>0} \times \Z_{>0}$.
We represent a diagram $D$ by drawing a cell in row $i$ 
column $j$ for each $(i, j) \in D$.
We use the matrix coordinates: Row $1$ is the topmost row
and column $1$ is the leftmost column. 
A \definition{weak composition}
is an infinite sequence
of $\Z_{\geq 0}$ with finitely many
positive entries. 
If $\alpha$ is a weak composition,
we use $\alpha_i$ to denote its $i\textsuperscript{th}$ entry. 
We write $\alpha$ as 
$(\alpha_1, \cdots, \alpha_n)$ where $\alpha_n$ is the last positive entry in $\alpha$.
The row weight (resp. column weight) of a diagram $D$ is a weak composition where the $i\textsuperscript{th}$
entry is the number of cells in row $i$ (resp. column $i$) of $D$.
We denote the row weight 
of a diagram $D$ by $\wt(D)$.

Pipedreams of size $n$
are in bijection with diagrams contained in $\{(i,j): 1 \leq i \leq n-1, 1 \leq j \leq n - i \}$.
Under this identification, 
$\PD(w)$ is a set of diagrams.
The ladder move is a move on diagrams 
and our algorithm is applying
ladder moves to diagrams.

\begin{exa}
We repeat Example~\ref{E: Algorithm} 
under our new convention:$$
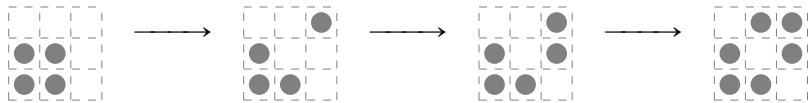
\begin{tikzpicture}[x=1em,y=1em,thick,color = blue]
    \draw[step=1,gray,ultra thin,dashed] (0,0) grid (3,3);
    \filldraw [gray] (1.5,0.5) circle (3.5pt);
    \filldraw [gray] (1.5,1.5) circle (3.5pt);
    \filldraw [gray] (0.5,1.5) circle (3.5pt);
    \filldraw [gray] (0.5,0.5) circle (3.5pt);
\end{tikzpicture}
\quad \raisebox{0.8cm}{$\xrightarrow{\quad\quad}$}\quad
\begin{tikzpicture}[x=1em,y=1em,thick,color = blue]
    \draw[step=1,gray,ultra thin,dashed] (0,0) grid (3,3);
    \filldraw [gray] (1.5,0.5) circle (3.5pt);
    \filldraw [gray] (2.5,2.5) circle (3.5pt);
    \filldraw [gray] (0.5,1.5) circle (3.5pt);
    \filldraw [gray] (0.5,0.5) circle (3.5pt);
\end{tikzpicture}
\quad \raisebox{0.8cm}{$\xrightarrow{\quad\quad}$}\quad
\begin{tikzpicture}[x=1em,y=1em,thick,color = blue]
    \draw[step=1,gray,ultra thin,dashed] (0,0) grid (3,3);
    \filldraw [gray] (1.5,0.5) circle (3.5pt);
    \filldraw [gray] (2.5,2.5) circle (3.5pt);
    \filldraw [gray] (2.5,1.5) circle (3.5pt);
    \filldraw [gray] (0.5,1.5) circle (3.5pt);
    \filldraw [gray] (0.5,0.5) circle (3.5pt);
\end{tikzpicture}
\quad \raisebox{0.8cm}{$\xrightarrow{\quad\quad}$}\quad
\begin{tikzpicture}[x=1em,y=1em,thick,color = blue]
\draw[step=1,gray,ultra thin,dashed] (0,0) grid (3,3);
    \filldraw [gray] (1.5,0.5) circle (3.5pt);
    \filldraw [gray] (1.5,2.5) circle (3.5pt);
    \filldraw [gray] (2.5,2.5) circle (3.5pt);
    \filldraw [gray] (2.5,1.5) circle (3.5pt);
    \filldraw [gray] (0.5,1.5) circle (3.5pt);
    \filldraw [gray] (0.5,0.5) circle (3.5pt);
\end{tikzpicture}
$$
The last diagram is $\widehat{P}(w)$
when $w$ has one-line notation $14523$.
Its row weight and column weight are
both $(2,2,2)$.
\end{exa}

\subsection{Snow diagrams
and $\rajcode$}

For any diagrams $D$,
Pan and Yu defined $\dark(D) \subseteq D$
which can be computed as follows:
Scan through $D$ from bottom to top.
For each row $r$, if there exists $(r,c) \in D$
such that currently there is no cells in column $c$
of $\dark(D)$,
we find the largest such $c$ and 
put $(r,c)$ in $\dark(D)$. 
Cells in $\dark(D)$ of $D$
are called \definition{dark clouds} of $D$.

\begin{exa}
The following is a diagram $D$
and $\dark(D)$
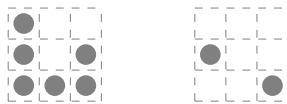
$$
\begin{tikzpicture}[x=1em,y=1em,thick,color = blue]
\draw[step=1,gray,ultra thin,dashed] (0,0) grid (3,3);
\filldraw [gray] (0.5,2.5) circle (3.5pt);
\filldraw [gray] (0.5,1.5) circle (3.5pt);
\filldraw [gray] (2.5,1.5) circle (3.5pt);
\filldraw [gray] (0.5,0.5) circle (3.5pt);
\filldraw [gray] (1.5,0.5) circle (3.5pt);
\filldraw [gray] (2.5,0.5) circle (3.5pt);
\end{tikzpicture} 
\quad\quad\quad
\begin{tikzpicture}[x=1em,y=1em,thick,color = blue]
\draw[step=1,gray,ultra thin,dashed] (0,0) grid (3,3);
\filldraw [gray] (0.5,1.5) circle (3.5pt);
\filldraw [gray] (2.5,0.5) circle (3.5pt);
\end{tikzpicture} 
$$
\end{exa}

There is an alternative characterization of
$\dark(D)$.
\begin{prop}
\label{P: Dark}
The diagram $\dark(D)$ is the unique subset of $D$
such that
\begin{itemize}
\item There is at most one cell in each row 
or column of $D$.
\item For any $(i,j) \in D$,
there is $(i',j) \in \dark(D)$ with $i' > i$
or there is $(i,j') \in \dark(D)$ with $j' > j$.
\end{itemize}
\end{prop}
\begin{proof}
By Remark 3.4 of~\cite{PY},
$\dark(D)$ satisfies the two conditions. 
The uniqueness is trivial. 
\end{proof}

The \definition{Rothe diagram} of $w$,
denoted as $\Rothe(w)$,
is the following diagram:
$$
\{(i,w(j)): i < j, w(i) > w(j)\}.
$$
For $w \in S_n$, the first $n$ 
numbers in $\wt(\Rothe(w))$
form $\invcode(w)$.
Let $\overleftarrow{\Rothe(w)}$
be the diagram obtained by left-justifying
all cells in $\Rothe(w)$.
This is the diagram in $\PD(w)$ 
that our algorithm starts with.

\begin{exa}
Take $w \in S_7$ with one-line notation $4617352$. The following are $\Rothe(w)$ and $\overleftarrow{\Rothe(w)}$.
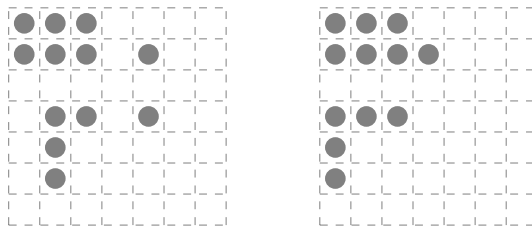
$$
\begin{tikzpicture}[x=1em,y=1em,thick,color = blue]
    \draw[step=1,gray,ultra thin,dashed] (0,0) grid (7,7);
    \filldraw [gray] (0.5,6.5) circle (3.5pt);
    \filldraw [gray] (0.5,5.5) circle (3.5pt);
    \filldraw [gray] (1.5,6.5) circle (3.5pt);
    \filldraw [gray] (1.5,5.5) circle (3.5pt);
    \filldraw [gray] (1.5,3.5) circle (3.5pt);
    \filldraw [gray] (1.5,2.5) circle (3.5pt);
    \filldraw [gray] (1.5,1.5) circle (3.5pt);
    \filldraw [gray] (2.5,6.5) circle (3.5pt);
    \filldraw [gray] (2.5,5.5) circle (3.5pt);
    \filldraw [gray] (2.5,3.5) circle (3.5pt);
    \filldraw [gray] (4.5,5.5) circle (3.5pt);
    \filldraw [gray] (4.5,3.5) circle (3.5pt);
\end{tikzpicture}
\quad\quad\quad
\begin{tikzpicture}[x=1em,y=1em,thick,color = blue]
    \draw[step=1,gray,ultra thin,dashed] (0,0) grid (7,7);
    \filldraw [gray] (0.5,6.5) circle (3.5pt);
    \filldraw [gray] (0.5,5.5) circle (3.5pt);
    \filldraw [gray] (1.5,6.5) circle (3.5pt);
    \filldraw [gray] (1.5,5.5) circle (3.5pt);
    \filldraw [gray] (0.5,3.5) circle (3.5pt);
    \filldraw [gray] (0.5,2.5) circle (3.5pt);
    \filldraw [gray] (0.5,1.5) circle (3.5pt);
    \filldraw [gray] (2.5,6.5) circle (3.5pt);
    \filldraw [gray] (2.5,5.5) circle (3.5pt);
    \filldraw [gray] (1.5,3.5) circle (3.5pt);
    \filldraw [gray] (3.5,5.5) circle (3.5pt);
    \filldraw [gray] (2.5,3.5) circle (3.5pt);
\end{tikzpicture}
$$
\end{exa}

For each $w \in S_n$,
Pechenik, Speyer and Weigandt
defined the weak composition 
$\rajcode(w)$
using increasing subsequences of $w$.
In this paper, 
we use a diagrammatic definition
of Pan and Yu~\cite{PY}

\begin{defn}[\cite{PY}]
\label{D: rajcode}  
Take $w \in S_n$ and find $\dark(\Rothe(w))$.
For each cell in $\dark(\Rothe(w))$, 
we fill all the empty cells above it
in $\Rothe(w)$. 
The resulting diagram is the \definition{snow diagram} of $w$.
Define $\rajcode(w)$ as the row weight
of the snow diagram of $w$.
\end{defn}

\begin{exa}
\label{E: rajcode}
Take $w \in S_7$ with one-line notation $4617352$. The following is its snow diagram. For clarity, we represent dark clouds
by a black circle and use
\textcolor{blue}{$\Asterisk$} to denote
the added cells.
$$
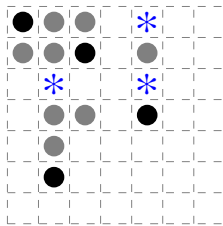
\begin{tikzpicture}[x=1em,y=1em,thick,color = blue]
    \draw[step=1,gray,ultra thin,dashed] (0,0) grid (7,7);
    \filldraw [black] (0.5,6.5) circle (3.5pt);
    \filldraw [gray] (0.5,5.5) circle (3.5pt);
    \filldraw [gray] (1.5,6.5) circle (3.5pt);
    \filldraw [gray] (1.5,5.5) circle (3.5pt);
    \filldraw [gray] (1.5,3.5) circle (3.5pt);
    \filldraw [gray] (1.5,2.5) circle (3.5pt);
    \filldraw [black] (1.5,1.5) circle (3.5pt);
    \filldraw [gray] (2.5,6.5) circle (3.5pt);
    \filldraw [black] (2.5,5.5) circle (3.5pt);
    \filldraw [gray] (2.5,3.5) circle (3.5pt);
    \filldraw [gray] (4.5,5.5) circle (3.5pt);
    \filldraw [black] (4.5,3.5) circle (3.5pt);
    \node[color=blue] at (1.5,4.5) {$\Asterisk$};
    \node[color=blue] at (4.5,4.5) {$\Asterisk$};
    \node[color=blue] at (4.5,6.5) {$\Asterisk$};
\end{tikzpicture}
$$
Thus, $\rajcode(w) = (4,4,2,3,1,1)$.
\end{exa}

It is well-known that $\Rothe(w)$
and $\Rothe(w^{-1})$ 
are conjugations of each other.
By Proposition~\ref{P: Dark},
$\dark(\Rothe(w))$ and $\dark(\Rothe(w^{-1}))$
are conjugations of each other.
Thus, we define the 
\definition{left snow diagram}
of $w$ as the diagram where
we fill empty spots
on the left of each dark cloud
in $\Rothe(w)$.
Its column weight will be
the same as the row weight
of the snow diagram of $w^{-1}$,
which is $\rajcode(w^{-1})$.

\begin{exa}
Keep the same $w$ as in Example~\ref{E: rajcode}.
Its left snow diagram is

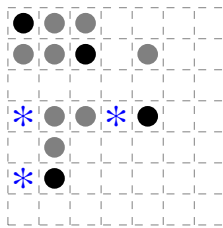
\begin{align*}
    \begin{tikzpicture}[x=1em,y=1em,thick,color = blue]
        \draw[step=1,gray,ultra thin,dashed] (0,0) grid (7,7);
        \filldraw [black] (0.5,6.5) circle (3.5pt);
        \filldraw [gray] (0.5,5.5) circle (3.5pt);
        \filldraw [gray] (1.5,6.5) circle (3.5pt);
        \filldraw [gray] (1.5,5.5) circle (3.5pt);
        \filldraw [gray] (1.5,3.5) circle (3.5pt);
        \filldraw [gray] (1.5,2.5) circle (3.5pt);
        \filldraw [black] (1.5,1.5) circle (3.5pt);
        \filldraw [gray] (2.5,6.5) circle (3.5pt);
        \filldraw [black] (2.5,5.5) circle (3.5pt);
        \filldraw [gray] (2.5,3.5) circle (3.5pt);
        \filldraw [gray] (4.5,5.5) circle (3.5pt);
        \filldraw [black] (4.5,3.5) circle (3.5pt);
        \node[color=blue] at (0.5,3.5) {$\Asterisk$};
        \node[color=blue] at (0.5,1.5) {$\Asterisk$};
        \node[color=blue] at (3.5,3.5) {$\Asterisk$};
    \end{tikzpicture}
\end{align*}
Thus, $\rajcode(w^{-1}) = (4,5,3,1,2)$.
\end{exa}

\section{Various recursions}
\label{S: Recursion}
We describe a recursive way to construct $\Rothe(w)$ 
and $\dark(\Rothe(w))$.
Then we obtain recursive formulas for 
$\rajcode(w)$ and $\rajcode(w^{-1})$.
Notice that $\invcode(\cdot)$
is a bijection from $S_n$ to 
weak compositions $(\alpha_1, \alpha_2, \cdots)$
where $\alpha_i \leq n - i$ for $i \in [n-1]$
and $\alpha_n = \alpha_{n+1} = \cdots = 0$.
We identify $w \in S_n$ with 
$(a, u) \in \{0, 1, \cdots, n-1\} \times S_{n-1}$ where
$a = \invcode(w)_1$ and $u$ 
is the unique permutation in $S_{n-1}$ with 
$\invcode(u) = (\invcode(w)_2,\invcode(w)_3, \cdots)$.
We simply write $w = (a,u)$.
Then we may recursively 
construct $\Rothe(w)$ as follows. 
Start from $\Rothe(u)$.
Shift all cells downward by $1$.
Then shift all cells in columns $a+1, a+2, \cdots$ to the right by $1$.
Finally, put cells at $(1,1), \cdots, (1,a)$.
The resulting diagram is $\Rothe(w)$.

Similarly, 
to construct $\dark(\Rothe(w))$,
we can start from $\dark(\Rothe(u))$.
Shift all cells downward by $1$.
Then shift all cells in columns $a+1, a+2, \cdots$ to the right by $1$.
Finally, find the largest $c \in [a]$ 
such that $\dark(\Rothe(u))$
has no cells in column $c$.
Put $(1,c)$ into $\dark(\Rothe(u))$.

\begin{exa}
\label{E: recursive Rothe}
    Keep $w \in S_7$ with one-line notation $4617352$.
    We have $w = (a,u)$ where $a = 3$ 
    and $u \in S_6$ has one-line notation $516342$.
    We depict how $\Rothe(u)$ and $\Rothe(w)$ as follows.
    The dark cells form $\dark(\Rothe(u))$
    and $\dark(\Rothe(w))$
    respectively. 
    \begin{align*}
        \begin{tikzpicture}[x=1em,y=1em,thick,color = blue]
            \draw[step=1,gray,ultra thin,dashed] (0,0) grid (6,6);
            \filldraw [gray] (0.5,5.5) circle (3.5pt);
            \filldraw [gray] (1.5,5.5) circle (3.5pt);
            \filldraw [gray] (1.5,3.5) circle (3.5pt);
            \filldraw [gray] (1.5,2.5) circle (3.5pt);
            \filldraw [black] (1.5,1.5) circle (3.5pt);
            \filldraw [black] (2.5,5.5) circle (3.5pt);
            \filldraw [gray] (2.5,3.5) circle (3.5pt);
            \filldraw [gray] (3.5,5.5) circle (3.5pt);
            \filldraw [black] (3.5,3.5) circle (3.5pt);
        \end{tikzpicture}
        &\quad \raisebox{1.9cm}{$\xrightarrow{\quad}$}\quad
        \begin{tikzpicture}[x=1em,y=1em,thick,color = blue]
            \draw[step=1,gray,ultra thin,dashed] (0,0) grid (7,7);
            \filldraw [black] (0.5,6.5) circle (3.5pt);
            \filldraw [gray] (1.5,6.5) circle (3.5pt);
            \filldraw [gray] (2.5,6.5) circle (3.5pt);
            \filldraw [gray] (0.5,5.5) circle (3.5pt);
            \filldraw [gray] (1.5,5.5) circle (3.5pt);
            \filldraw [gray] (1.5,3.5) circle (3.5pt);
            \filldraw [gray] (1.5,2.5) circle (3.5pt);
            \filldraw [black] (1.5,1.5) circle (3.5pt);
            \filldraw [black] (2.5,5.5) circle (3.5pt);
            \filldraw [gray] (2.5,3.5) circle (3.5pt);
            \filldraw [gray] (4.5,5.5) circle (3.5pt);
            \filldraw [black] (4.5,3.5) circle (3.5pt);
        \end{tikzpicture}
    \end{align*}
\end{exa}

Consequently, we may compute $\rajcode(w)$ and $\rajcode(w^{-1})$ recursively.
Let $d_c(u)$ be the number of cells
in $\dark(\Rothe(u))$ that are strictly to 
the right of column $c$.

\begin{prop}
\label{P: Recursive rajcode}   
Take $w = (a,u) \in S_n$.
\begin{itemize}
\item We can get $\rajcode(w)$ by
prepending $a + d_a(u)$ to $\rajcode(u)$.
\item To obtain $\rajcode(w^{-1})$,
we just insert $d_a(u)$ between the $a\textsuperscript{th}$
and $(a+1)\textsuperscript{th}$ entries of $\rajcode(u^{-1})$.
Then increase the first $a$ entries by $1$.
\end{itemize}
Consequently, $\reg(w) - \reg(u) = d_a(u)$.
\end{prop}

\begin{proof}
Follows directly from the recursive constructions
of $\Rothe(w)$ and $\dark(\Rothe(w))$.
\end{proof}

\begin{exa}
Keep $w = (a,u)$ in Example~\ref{E: recursive Rothe}.
We show how the snow diagram and left snow diagram 
of $w$ differ from those of $u$:
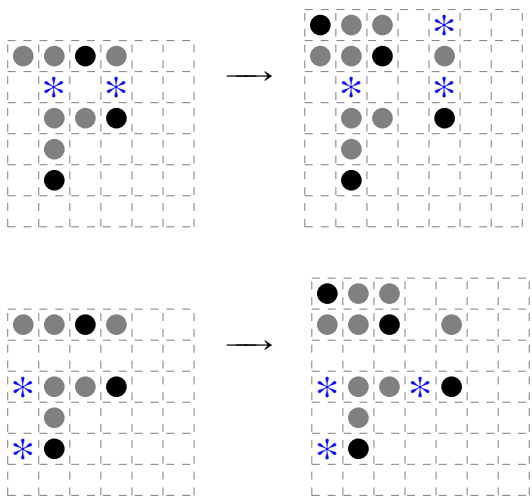
    \begin{align*}
        \begin{tikzpicture}[x=1em,y=1em,thick,color = blue]
            \draw[step=1,gray,ultra thin,dashed] (0,0) grid (6,6);
            \filldraw [gray] (0.5,5.5) circle (3.5pt);
            \filldraw [gray] (1.5,5.5) circle (3.5pt);
            \filldraw [gray] (1.5,3.5) circle (3.5pt);
            \filldraw [gray] (1.5,2.5) circle (3.5pt);
            \filldraw [black] (1.5,1.5) circle (3.5pt);
            \filldraw [black] (2.5,5.5) circle (3.5pt);
            \filldraw [gray] (2.5,3.5) circle (3.5pt);
            \filldraw [gray] (3.5,5.5) circle (3.5pt);
            \filldraw [black] (3.5,3.5) circle (3.5pt);
            \node[color=blue] at (1.5,4.5) {$\Asterisk$};
            \node[color=blue] at (3.5,4.5) {$\Asterisk$};
        \end{tikzpicture}
        &\quad \raisebox{1.9cm}{$\xrightarrow{\quad}$}\quad
        \begin{tikzpicture}[x=1em,y=1em,thick,color = blue]
            \draw[step=1,gray,ultra thin,dashed] (0,0) grid (7,7);
            \filldraw [black] (0.5,6.5) circle (3.5pt);
            \filldraw [gray] (1.5,6.5) circle (3.5pt);
            \filldraw [gray] (2.5,6.5) circle (3.5pt);
            \filldraw [gray] (0.5,5.5) circle (3.5pt);
            \filldraw [gray] (1.5,5.5) circle (3.5pt);
            \filldraw [gray] (1.5,3.5) circle (3.5pt);
            \filldraw [gray] (1.5,2.5) circle (3.5pt);
            \filldraw [black] (1.5,1.5) circle (3.5pt);
            \filldraw [black] (2.5,5.5) circle (3.5pt);
            \filldraw [gray] (2.5,3.5) circle (3.5pt);
            \filldraw [gray] (4.5,5.5) circle (3.5pt);
            \filldraw [black] (4.5,3.5) circle (3.5pt);
            \node[color=blue] at (1.5,4.5) {$\Asterisk$};
            \node[color=blue] at (4.5,4.5) {$\Asterisk$};
            \node[color=blue] at (4.5,6.5) {$\Asterisk$};
        \end{tikzpicture}
    \end{align*}
    \begin{align*}
        \begin{tikzpicture}[x=1em,y=1em,thick,color = blue]
            \draw[step=1,gray,ultra thin,dashed] (0,0) grid (6,6);
            \filldraw [gray] (0.5,5.5) circle (3.5pt);
            \filldraw [gray] (1.5,5.5) circle (3.5pt);
            \filldraw [gray] (1.5,3.5) circle (3.5pt);
            \filldraw [gray] (1.5,2.5) circle (3.5pt);
            \filldraw [black] (1.5,1.5) circle (3.5pt);
            \filldraw [black] (2.5,5.5) circle (3.5pt);
            \filldraw [gray] (2.5,3.5) circle (3.5pt);
            \filldraw [gray] (3.5,5.5) circle (3.5pt);
            \filldraw [black] (3.5,3.5) circle (3.5pt);
            \node[color=blue] at (0.5,3.5) {$\Asterisk$};
            \node[color=blue] at (0.5,1.5) {$\Asterisk$};
        \end{tikzpicture}
        &\quad \raisebox{1.9cm}{$\xrightarrow{\quad}$}\quad
        \begin{tikzpicture}[x=1em,y=1em,thick,color = blue]
            \draw[step=1,gray,ultra thin,dashed] (0,0) grid (7,7);
            \filldraw [black] (0.5,6.5) circle (3.5pt);
            \filldraw [gray] (1.5,6.5) circle (3.5pt);
            \filldraw [gray] (2.5,6.5) circle (3.5pt);
            \filldraw [gray] (0.5,5.5) circle (3.5pt);
            \filldraw [gray] (1.5,5.5) circle (3.5pt);
            \filldraw [gray] (1.5,3.5) circle (3.5pt);
            \filldraw [gray] (1.5,2.5) circle (3.5pt);
            \filldraw [black] (1.5,1.5) circle (3.5pt);
            \filldraw [black] (2.5,5.5) circle (3.5pt);
            \filldraw [gray] (2.5,3.5) circle (3.5pt);
            \filldraw [gray] (4.5,5.5) circle (3.5pt);
            \filldraw [black] (4.5,3.5) circle (3.5pt);
            \node[color=blue] at (0.5,3.5) {$\Asterisk$};
            \node[color=blue] at (0.5,1.5) {$\Asterisk$};
            \node[color=blue] at (3.5,3.5) {$\Asterisk$};
        \end{tikzpicture}
    \end{align*}
We have $d_a(u) = 1$.
We obtain $\rajcode(w) = (4, 4, 2, 3, 1, 1)$
by prepending $a + d_a(u) = 4$
to $\rajcode(u) = (4,2,3,1,1)$.
We obtain $\rajcode(w^{-1}) = (4,5,3,1,2)$
by inserting $d_a(u)$ after the $a\textsuperscript{th}$ entry of $\rajcode(u^{-1}) = (3,4,2,2)$ and then increase the first $a$
entries by $1$.
\end{exa}

Notice that when $w = (a,u)$,
$\invcode(w)$ can be obtained
by prepending the number $a$ to $\invcode(u)$.
Thus, we also have a recursive formula
for $\IR(w)$.
For a monomial $m$,
let $\overrightarrow{m}$ be the monomial
obtained by turning each $x_i$ in $m$
into $x_{i+1}$.

\begin{prop}
~\label{P: recursion}
Take $w = (a,u) \in S_n$.
Let $$(M_0, \cdots, M_{\reg(w)}) = \IR(w), (m_0, \cdots, m_{\reg(u)}) = \IR(u).$$
Then $\reg(w) = \reg(u) + d_a(u)$ and 
$$
M_j = \begin{cases}
x_1^a \overrightarrow{m_j} & \text{if $j = 0, 1, \cdots, \reg(u)$,} \\
x_1^{a + j - \reg(u)} \times \overrightarrow{m_{\reg(u)}} & \text{if $j = \reg(u) + 1, \cdots, \reg(w)$.}
\end{cases}$$
\end{prop}
\begin{proof}
Follows directly from the recursive formula of
$\rajcode(w)$ and the definition of $\IR(\cdot)$.
\end{proof}

\begin{exa}
    Keep $w = (a,u)$ in Example~\ref{E: recursive Rothe}. 
    We have $\reg(u) = 2$ and $\reg(w) = \reg(u) + d_a(u) = 3$.
    Since 
    \begin{align*}
        \IR(u) = (x^{(4,0,3,1,1)}, x^{(4,1,3,1,1)}, x^{(4,2,3,1,1)}),
    \end{align*}
    we have
    \begin{align*}
        \IR(w) = (x^{(3,4,0,3,1,1)}, x^{(3,4,1,3,1,1)}, x^{(3,4,2,3,1,1)},  x^{(4,4,2,3,1,1)})
    \end{align*}

\end{exa}

\section{Proof of main theorems}
\label{S: Proofs}

To prove our main theorems,
we need to introduce a 
new permutation statistic.
\begin{defn}
\label{D: movecode} 
For $w \in S_n$,
its \definition{movecode}, 
denoted as $\movecode(w)$,
is a weak composition where
$\movecode(w)_{i}$ is the number of cells in column $i$ of $\Rothe(w)$ with no dark clouds strictly to its right.
\end{defn}

\begin{exa}
    Take $w \in S_7$ with one-line notation $4617352$. The following is $\Rothe(w)$, where the black cells are dark clouds and blue cells are non-dark cloud cells without dark clouds to their right.
    $$
    
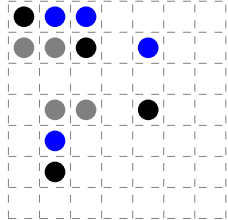
\begin{tikzpicture}[x=1em,y=1em,thick,color = blue]
        \draw[step=1,gray,ultra thin,dashed] (0,0) grid (7,7);
        \filldraw [black] (0.5,6.5) circle (3.5pt);
        \filldraw [gray] (0.5,5.5) circle (3.5pt);
        \filldraw [blue] (1.5,6.5) circle (3.5pt);
        \filldraw [gray] (1.5,5.5) circle (3.5pt);
        \filldraw [gray] (1.5,3.5) circle (3.5pt);
        \filldraw [blue] (1.5,2.5) circle (3.5pt);
        \filldraw [black] (1.5,1.5) circle (3.5pt);
        \filldraw [blue] (2.5,6.5) circle (3.5pt);
        \filldraw [black] (2.5,5.5) circle (3.5pt);
        \filldraw [gray] (2.5,3.5) circle (3.5pt);
        \filldraw [blue] (4.5,5.5) circle (3.5pt);
        \filldraw [black] (4.5,3.5) circle (3.5pt);
    \end{tikzpicture}
    $$
    Then $\movecode(w)$ is the number of black and blue cells in each column, which is $(1,3,2,0,2)$.
\end{exa}

We have the following observation regarding
this permutation statistic.
\begin{prop}
\label{P: movecode and d}
Take $w \in S_n$
and $c \in [n]$.
Then 
$$\rajcode(w^{-1})_{c+1} - \max(\movecode(w)_{c+1} - 1, 0) = d_c(u) 
= \rajcode(w^{-1})_{c} -\movecode(w)_{c}.$$
\end{prop}

\begin{proof}
We refer to cells in $\dark(\Rothe(w))$
as dark clouds. 
Consider the left snow diagram
of $w$. 
In the diagram, 
there are four types of cells. 
\begin{itemize}
\item Type 1: Dark clouds
\item Type 2: Cells that do not belong 
to $\Rothe(w)$.
\item Type 3: Cells in $\Rothe(w)$ with a dark
cloud in its row on its right.
\item Type 4: Cells in $\Rothe(w)$ that is not
a dark cloud and has no dark
cloud in its row on its right.
\end{itemize}

The number of 
type $1$, $2$ and $4$ cells in column $c+1$
is $d_c(w)$.
The number of all cells in column $c+1$
is $\rajcode(w^{-1})_{c+1}$.
The number of type $3$ cells in column $c+1$
is $\max(\movecode(w)_c - 1, 0)$,
so we have the first equation. 

The number of 
type $2$ and $3$ cells in column $c$
is $d_c(w)$.
The number of all cells in column $c$
is $\rajcode(w^{-1})_{c}$.
The number of type $1$ and $4$ 
cells in column $c$
is $\movecode(w)_c$,
so we have the second equation. 
\end{proof}

The main application of $\movecode(w)$
is to characterize the number of cells
moved when our algorithm processes each column.

\begin{prop}
\label{P: Movecode}  
Take $v = (a,w) \in S_n$.
During the last iteration of the algorithm 
that computes $\widehat{P}(v)$,
the number of cells moved in column $c$ is $\movecode(w)_c$ if $c > a$ and $0$ otherwise.
\end{prop}

\begin{exa}
    Keep $v \in S_7$ with one-line notation $4617352$.
    We have $v = (a, w)$
    where $a = 3$ and $w \in S_6$
    has one-line notation 
    $516342$.
    We have $\movecode(w) = (0,2,1,2)$. During the last iteration of the algorithm, 
    the bar is right above row $1$.
    The algorithm moves $0$ cells in column $c > 4$,
    since $\movecode(w)_c = 0$.
    The algorithm moves $2$ cells in column $4$ since $4 > a$ and $\movecode(w)_4 = 2$.
    It moves $0$ cells in column $3, 2,$ and $1$
    since $1, 2, 3 \leq a$.
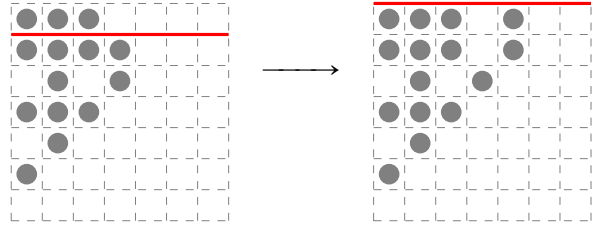
    $$
    \begin{tikzpicture}[x=1em,y=1em,thick,color = blue]
        \draw[step=1,gray,ultra thin,dashed] (0,0) grid (7,7);
        \draw[color=red, very thick] (0,6)--(7,6);
        \filldraw [gray] (0.5,6.5) circle (3.5pt);
        \filldraw [gray] (0.5,5.5) circle (3.5pt);
        \filldraw [gray] (0.5,3.5) circle (3.5pt);
        \filldraw [gray] (0.5,1.5) circle (3.5pt);
        \filldraw [gray] (1.5,6.5) circle (3.5pt);
        \filldraw [gray] (1.5,5.5) circle (3.5pt);
        \filldraw [gray] (1.5,4.5) circle (3.5pt);
        \filldraw [gray] (1.5,3.5) circle (3.5pt);
        \filldraw [gray] (1.5,2.5) circle (3.5pt);
        \filldraw [gray] (2.5,6.5) circle (3.5pt);
        \filldraw [gray] (2.5,5.5) circle (3.5pt);
        \filldraw [gray] (2.5,3.5) circle (3.5pt);
        \filldraw [gray] (3.5,5.5) circle (3.5pt);
        \filldraw [gray] (3.5,4.5) circle (3.5pt);
    \end{tikzpicture}
    \quad \raisebox{1.9cm}{$\xrightarrow{\quad\quad}$}\quad
    \begin{tikzpicture}[x=1em,y=1em,thick,color = blue]
        \draw[step=1,gray,ultra thin,dashed] (0,0) grid (7,7);
        \draw[color=red, very thick] (0,7)--(7,7);
        \filldraw [gray] (0.5,6.5) circle (3.5pt);
        \filldraw [gray] (0.5,5.5) circle (3.5pt);
        \filldraw [gray] (0.5,3.5) circle (3.5pt);
        \filldraw [gray] (0.5,1.5) circle (3.5pt);
        \filldraw [gray] (1.5,6.5) circle (3.5pt);
        \filldraw [gray] (1.5,5.5) circle (3.5pt);
        \filldraw [gray] (1.5,4.5) circle (3.5pt);
        \filldraw [gray] (1.5,3.5) circle (3.5pt);
        \filldraw [gray] (1.5,2.5) circle (3.5pt);
        \filldraw [gray] (2.5,6.5) circle (3.5pt);
        \filldraw [gray] (2.5,5.5) circle (3.5pt);
        \filldraw [gray] (2.5,3.5) circle (3.5pt);
        \filldraw [gray] (3.5,4.5) circle (3.5pt);
        \filldraw [gray] (4.5,6.5) circle (3.5pt);
        \filldraw [gray] (4.5,5.5) circle (3.5pt);
    \end{tikzpicture}
    $$
\end{exa}

We prove this proposition in \S\ref{S: Proofs of Prop and Cor}. 
Our proof requires a few technical 
lemmas which also lead to the following
result:

\begin{cor}
\label{C: Row +1}
Consider the iteration when 
the bar is right above row $i$
in our algorithm.
Let $D_1$ (resp. $D_2$) 
be the diagram before (resp. after)
processing one column.
If the algorithm makes a move
in this column, 
then $\wt(D_2)$ is obtained
from increasing $i\textsuperscript{th}$
entry of $\wt(D_1)$ by $1$.
\end{cor}

Using Proposition~\ref{P: Movecode}
and Corollary~\ref{C: Row +1},
we can prove our main results. 
We start with 
Theorem~\ref{T: Intermidiate row}.

\begin{proof}[Proof of Theorem~\ref{T: Intermidiate row}]
      We induct on $n$. The base case ($n = 1$) is trivial. Let $w = (a,u) \in S_{n}$ with $n > 1$. 
By our inductive hypothesis,
the algorithm made $\reg(u)$ K-ladder
moves before the last iteration.
By Proposition~\ref{P: Movecode}, 
in the last iteration of the algorithm,
it makes a K-ladder move in column $c$
if and only if $c > a$ and $\movecode(u)_c > 0$.
This is exactly the number $d_a(u)$,
which equals $\reg(w) - \reg(u)$
by Proposition~\ref{P: Recursive rajcode}.
Thus, the algorithm to compute $\widehat{P}(w)$
makes $\reg(w)$ K-ladder moves in total.

      Let
      $$\IR(w) = (M_0, \cdots, M_{\reg(w)}), \IR(u) = (m_0, \cdots, m_{\reg(u)}).$$
      By Proposition~\ref{P: recursion}, for $i = 0, \cdots, \reg(u)$,
      we have $M_i = x_1^a \overrightarrow{m_i}$.
      When the algorithm
      makes the $i\textsuperscript{th}$ K-ladder
      move, the bar has not reached row $1$.
      Before the bar reaches row $1$,
      the algorithm ignores the first row 
      of the diagram, which has $a$ cells, 
      and behaves as if computing
      $\widehat{P}(u)$.
      Thus, the statement holds for $i  = 0, 1, \cdots, \reg(u)$
      by our inductive hypothesis.
    
      For $i = \reg(u) + 1, \cdots, \reg(w)$, 
      the $i\textsuperscript{th}$ K-ladder move
      happens when the bar is above row $1$.
      Let $D$ be the diagram right after the $(i-1)\textsuperscript{th}$ K-ladder move and $D'$ be the diagram right after the $i\textsuperscript{th}$ K-ladder move. By Corollary~\ref{C: Row +1}, $x^{\wt(D')} = x_1 \cdot x^{\wt(D)}$, which concludes the proof. 
\end{proof}

\begin{proof}[Proof of Theorem~\ref{T: Row col weight}]
    By Theorem~\ref{T: Intermidiate row},
    the row weight of $\widehat{P}(w)$ is $\rajcode(w)$.
    For the column weight, 
    we prove by induction on $n$.
The base case $n = 1$ is trivial.
Now assume $n > 1$ and $w = (a,u) \in S_n$.
Let $D$ be the diagram we have right before
the last iteration of the algorithm 
computing $\widehat{P}(w)$. 
It can be obtained by shifting $\widehat{P}(u)$
downward by $1$ and append $a$ left-justified cells
in the first row. 
By our inductive hypothesis, 
$\widehat{P}(u)$ has column weight $\rajcode(u^{-1})$.
Now take $c \in [n-1]$ and consider three cases:
\begin{itemize}
\item Suppose $c > a + 1$.
Consider the last iteration of the algorithm.
By Proposition~\ref{P: Movecode},
the algorithm makes $\movecode(u)_c$ 
(resp. $\movecode(u)_{c-1}$) moves in column $c$ (resp. $c-1$).
Thus, column $c$ loses $\max(\movecode(u)_c-1, 0)$
cells and then gain $\movecode(u)_{c-1}$ cells. 
By Proposition~\ref{P: movecode and d},
$\widehat{P}(w)$ has 
$$
\rajcode_c(u^{-1}) - \max(\movecode(u)_c-1, 0) 
+ \movecode(u)_{c-1} = \rajcode_{c-1}(u^{-1})
$$
cells in column $c$.
Finally, by 
Proposition~\ref{P: Recursive rajcode}, 
$\rajcode_{c-1}(u^{-1})$ is just $\rajcode_{c}(w^{-1})$.
\item Suppose $c = a + 1$.
By Proposition~\ref{P: Movecode},
the algorithm makes $\movecode(u)_c$ 
moves in column $c$,
and makes $0$ moves in column $c-1$
if it exists.
Thus, column $c$ loses $\max(\movecode(u)_c-1, 0)$
cells. 
By Proposition~\ref{P: movecode and d},
$\widehat{P}(w)$ has 
$$
\rajcode_c(u^{-1}) - \max(\movecode(u)_c-1, 0) = 
d_a(u)
$$
cells in column $c$.
Finally, by 
Proposition~\ref{P: Recursive rajcode}, 
$d_a(u)$ is just $\rajcode_{c}(w^{-1})$.
\item Suppose $c \in [a]$.
By Proposition~\ref{P: Movecode},
the algorithm makes $0$ 
moves in column $c$,
and makes $0$ moves in column $c-1$
if it exists.
Thus, $\widehat{P}(w)$ 
has $\rajcode(u^{-1})_c + 1$
cells in column $c$.
Finally, by 
Proposition~\ref{P: Recursive rajcode}, 
$\rajcode(u^{-1})_c + 1$ is just $\rajcode_{c}(w^{-1})$. \qedhere
\end{itemize}
\end{proof}

\section{Proof of Proposition~\ref{P: Movecode}
and Corollary~\ref{C: Row +1}}
\label{S: Proofs of Prop and Cor}

Following \S\ref{S: Recursion},
we derive a recursive way to compute
$\movecode(w)$.

\begin{lem}
    ~\label{L: constructmovecode}
    For $w \in S_n$, 
    we write $w= (a,u)$.
    Then $\movecode(w)$ can be determined starting from $\movecode(u)$. 
    First, insert a $0$ between 
    $\movecode(u)_a$ and $\movecode(u)_{a+1}$.
    Then start from the $a\textsuperscript{th}$ entry
    and increase each entry by $1$ from right to left. 
    Whenever we change a $0$ into a $1$, we stop immediately. 
    The resulting weak composition is $\movecode(w)$.
\end{lem}

\begin{proof}
Follows directly from the recursive constructions
of $\Rothe(w)$ and $\dark(\Rothe(w))$.
\end{proof}

\begin{exa}
    Take $w \in S_7$ with one-line notation $4617352$.
    We have $w = (3, u)$
    where $u \in S_6$
    has one-line notation 
    $516342$.
    We have $\movecode(u) = (0,2,1,2)$.
    Then we insert a $0$
    between $\movecode(u)_3$
    and $\movecode(u)_4$,
    obtaining $(0,2,1,0,2)$.
    We then increases entries
    by $1$ from right to left, 
    starting from the thrid entry.
    When we turn the $0$ in the first entry into $1$,
    we stop,
    obtaining $(1,3,2,0,2)$.
\end{exa}

Our proofs rely 
on a simple operator
on diagrams. 
We may break the algorithm
into a sequence of this operator.

\begin{defn}
\label{D: lift} 
We define the operator $L_{i,c}$
on diagrams. 
Take diagram $D$ and put 
a bar above row $i$ in $D$.
We ignore everything above the bar,
imagining row $i$ is the top-most row. 
Then we scan through cells in column 
$c$ from top to bottom. 
Whenever we see a cell at which
we can perform a ladder move,
we perform a regular ladder move. 
After going through this column, 
if we made a move, turn the last move
into a K-ladder move.
\end{defn}

With this notion, applying the algorithm on $w \in S_n$ can be rewritten as
\begin{align}
\label{EQ: algorithm}
\widehat{P}(w) = 
(L_{1, 1} \cdots L_{1, n-2})
\cdots
(L_{n-3, 1}L_{n-3, 2})(L_{n-2, 1})(\overleftarrow{\Rothe(w)})
\end{align}
In words, we iterate through
$i = n-2, \cdots, 2, 1$.
For each $i$,
we iterate through $c = n-1 - i, \cdots, 2, 1$ and apply $L_{i,c}$.

We start by observing a straightforward recursive 
property of this operator. 
\begin{rem}
\label{R: L recursion}
Fix $i,c \in \mathbb{Z}_{>0}$
and let $D$ be a diagram.
Suppose $(i,c) \notin D$ and $(i, c+1) \notin D$.
\begin{itemize}
\item Suppose $(i+1, c) \in D$ and $(i+1, c+1) \notin D$.
Let $D'$ be the diagram obtained by moving $(i+1,c)$
to $(i, c+1)$ in $D$.
If $L_{i+1, c}(D') \neq D'$,
we know $L_{i,c}(D) = L_{i+1, c}(D')$.
Otherwise, $L_{i,c}(D) = D' \sqcup \{(i+1,c)\}$.
Informally, in this case, 
$L_{i,c}$ behaves as if $L_{i+1,c}$
after the regular ladder move on $(i+1,c)$.
\item Suppose $(i+1, c) \in D$ and $(i+1, c+1) \in D$.
Then intuitively, 
$L_{i,c}$ behaves as if row $i+1$ is ignored:
Let $D'$ be obtained from $D$ by removing $(i+1, c)$
and $(i+1, c+1)$.
If $(i+1, c+1) \notin L_{i+1, c}(D')$,
$L_{i,c}(D) = L_{i+1, c}(D') \sqcup \{(i+1,c), (i+1, c+1)\}$.
Otherwise, 
$L_{i,c}(D) = L_{i+1, c}(D') \sqcup \{(i+1,c), (i, c+1)\}$.
\end{itemize}
\end{rem}

We are primarily interested in applying 
$L_{i,c}$ to a diagram
in the following case. 

\begin{defn}
We say the operator $L_{i,c}$ \definition{acts initially on}
$D$ if $D$ is fixed by $L_{i+1,c}$.
\end{defn}

Eventually, we will show all $L_{i,c}$
in our algorithm acts initially. 
We first derive a few properties 
when $L_{i,c}$ acts initially on $D$.

\begin{lem}
    \label{L: chainmove}
    Suppose $L_{i,c}$ acts initially on $D$ 
    and $L_{i,c}$ moves at least one cell.
    We let $(r_1, c), \cdots, (r_k, c)$ be the cells
    moved where $r_1 < \cdots < r_k$. 
    Let $r_0 = i$.
    Then we know the cell $(r_j, c)$ is moved to $(r_{j-1}, c+1)$ for $j \in [k]$. 
Thus, $\wt(L_{i,c}(D))$ is obtained from $\wt(D)$
    by adding $1$ to the $i\textsuperscript{th}$ entry.
\end{lem}
\begin{proof}
    If $L_{i,c}$ moves $(r_1,c)$ to $(r', c+1)$ 
    for some $r' > i$, then $L_{i+1,c}$ will also
    move $(r_1,c)$ to $(r', c+1)$.
    This contradicts our assumption
    that $L_{i,c}$ acts initially on $D$.
    Thus, $L_{i,c}$ moves $(r_1,c)$ to $(i, c+1)$. 
    
    For $j > 1$, when $(r_j, c)$ moves, $(r_{j-1},c)$ and $(r_{j-1},c+1)$ must both be empty since the cell in $(r_{j-1},c)$ just performed a ladder move. Therefore $(r_j, c)$ must be moved to $(r',c+1)$ for some $r' \geq r_{j-1}$. However, $r' > r_{j-1}$ contradicts the assumption that $L_{i,c}$ acts initially on $D$, so $r' = r_{j-1}$. 
\end{proof}

To better describe the effect of $L_{i,c}$
when it acts initially,
we introduce the following notion.

\begin{defn}
\label{D: initialsegment}  
The \definition{$(i,c)$-initial segment} 
of a diagram $D$ is the set of $(r,c)$
such that $(r',c) \in D$ for all $i \leq r' \leq r$.
\end{defn}

This notion characterizes the destination 
of cells moved by $L_{i,c}$ when it acts initially.

\begin{lem}
\label{L: intoinitial}
Suppose $L_{i,c}$ acts initially on $D$.
Then it moves cells 
to the $(i,c+1)$-initial segment of $L_{i,c}(D)$. 
\end{lem}

\begin{proof}
    Let $(r_1, c), (r_2, c), \dots, (r_k, c)$ where $r_1 < r_2< \cdots < r_k$ be the cells of $D$ moved by $L_{i,c}$. 
    Let $r_0 = i$.
    By Lemma~\ref{L: chainmove}, for $j \in [k]$,
    $(r_j, c)$ is moved to $(r_{j-1}, c+1)$.
    We show $(r_{j-1},c)$ is in the $(j, c+1)$-initial segment of $L_{i,c}(D)$ by induction on $j$. 
    For the base case, $(r_0, c+1) = (i, c+1)$ 
    is clearly in the $(j, c+1)$-initial segment of $L_{i,c}(D)$
    
    For $j > 1$.
    assume $(r_{j-2}, c + 1)$ is in the $(i, c+1)$-initial segment of $L_{i,c}(D)$.
    Since $(r_{j-1}, c)$ is moved to $(r_{j-2}, c + 1)$,
    we know $(r', c+1) \in L_{i,c}(D)$
    for any $r_{j-2} < r' < r_{j-1}$.
    Thus, $(r_{j-1}, c + 1)$ is in the $(i, c+1)$-initial segment of $L_{i,c}(D)$. 
\end{proof}

We can also use ``initial segment'' 
to characterize what cells can be moved
by $L_{i,c}$ when it acts initially.

\begin{lem}
    \label{L: movablecell}
     Suppose $L_{i,c}$ acts initially on $D$.
     If $(i,c) \in D$,
     then $D$ is fixed by $L_{i,c}$.
     Otherwise, a cell $(r,c) \in D$ is moved by $L_{i,c}$ if and only if it is in the $(i+1,c)$-initial segment of $D$ and $(r, c+1) \notin D$. 
\end{lem}
\begin{proof}
The lemma is immediate when $(i,c) \in D$.
Otherwise, let $(r_1, c), \cdots, (r_k, c) \in D$ 
be the cells moved by $L_{i,c}$
where $r_1 < \cdots < r_k$.
Let $r_0 = i$.
Clearly, $(r_j, c+1) \notin D$
for each $j \in [k]$.
We prove $(r_j,c)$
is in the $(i+1, c)$-initial segment of $D$
by induction. 
First, by Lemma~\ref{L: chainmove},
$(r_1, c)$ is moved 
to $(r_0, c+1)$,
so $(r', c) \in D$ for $r_0 = i < r' < r_1$.
In other words,
$(r_1, c)$ is in the $(i+1, c)$-initial segment of $D$.
For $j > 1$, by Lemma~\ref{L: chainmove},
$(r_j, c)$ is moved 
to $(r_{j-1}, c+1)$,
so $(r', c) \in D$ for $r_{j-1} < r' < r_j$.
The inductive step is finished
since $(r_{j-1}, c)$ is in the $(i+1, c)$-initial segment of $D$.

Now assume $(r, c)$ is a cell in the $(i+1, c)$-initial segment of $D$
and $(r, c+1) \notin D$.
Assume toward contradiction that
$(r,c)$ is not moved by $L_{i,c}$.
Take the smallest such $r$. 
Since $L_{i,c}$ moves $(r_j, c)$ to $(r_{j-1},c)$,
we know $(r', c+1) \in D$ for any $r_{j-1} < r' < r_j$.
Thus, we cannot have $r_{j-1} < r < r_j$
for $j \in [k]$.
Since $(r,c)$ is not moved, 
we know $r$ is not $r_1, \cdots, r_k$.
Thus, $r > r_k$.
By the minimality of $r$, 
$(r', c), (r', c+1) \in D$
for $r_k < r' < r$.
Thus, $L_{i,c}$ moves $(r_k, c)$, 
it can perform a ladder move at $(r,c)$. 
Contradiction.
\end{proof}

The following example is a demonstration of the previous two lemmas related to initial segments.

\begin{exa}
    Let $D$ be a diagram whose column $3$ and $4$ look like
    the picture on the left. 
    Notice that $D$ will be fixed by $L_{2,3}$.
    After applying $L_{1, 3}$, these two columns look like
    the picture on the right:
    $$
    
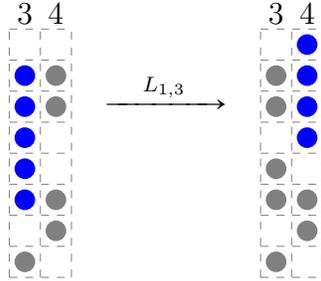
\begin{tikzpicture}[x=1em,y=1em,thick,color = blue]
        \draw[step=1,gray,ultra thin,dashed] (0,0) grid (2,8);
        \filldraw [gray] (1.5,6.5) circle (3.5pt);
        \filldraw [blue] (0.5,6.5) circle (3.5pt);
        \filldraw [gray] (1.5,5.5) circle (3.5pt); 
        \filldraw [blue] (0.5,5.5) circle (3.5pt); 
        \filldraw [blue] (0.5,4.5) circle (3.5pt); 
        \filldraw [blue] (0.5,3.5) circle (3.5pt); 
        \filldraw [blue] (0.5,2.5) circle (3.5pt); 
        \filldraw [gray] (1.5,2.5) circle (3.5pt);
        \filldraw [gray] (1.5,1.5) circle (3.5pt);
        \filldraw [gray] (0.5,0.5) circle (3.5pt);
        \node[color=black] at (0.5,8.5) {$3$};
        \node[color=black] at (1.5,8.5) {$4$};
    \end{tikzpicture}
    \quad \raisebox{2.2cm}{$\xrightarrow{\quad L_{1,3} \quad}$}\quad
    \begin{tikzpicture}[x=1em,y=1em,thick,color = blue]
        \draw[step=1,gray,ultra thin,dashed] (0,0) grid (2,8);
        \filldraw [blue] (1.5,7.5) circle (3.5pt); 
        \filldraw [blue] (1.5,6.5) circle (3.5pt);
        \filldraw [gray] (0.5,6.5) circle (3.5pt);
        \filldraw [blue] (1.5,5.5) circle (3.5pt); 
        \filldraw [gray] (0.5,5.5) circle (3.5pt); 
        \filldraw [blue] (1.5,4.5) circle (3.5pt);
        \filldraw [gray] (0.5,3.5) circle (3.5pt); 
        \filldraw [gray] (0.5,2.5) circle (3.5pt); 
        \filldraw [gray] (1.5,2.5) circle (3.5pt);
        \filldraw [gray] (1.5,1.5) circle (3.5pt);
        \filldraw [gray] (0.5,0.5) circle (3.5pt);
        \node[color=black] at (0.5,8.5) {$3$};
        \node[color=black] at (1.5,8.5) {$4$};
    \end{tikzpicture}
    $$
    We color the $(2,3)$-initial segment of $D$
    and $(1,4)$-initial segment of $L_{1,3}(D)$.
    Notice that $L_{1,3}$ move cells to the 
    $(1,4)$-initial segment of $L_{1,3}(D)$. Also notice that cells in column $3$ is moved if and only if it is in the $(2,3)$-initial segment of $D$ and has no cell on its right.
\end{exa}

We also have the ``converse statement'' of Lemma~\ref{L: movablecell}.
\begin{lemma}
\label{L: initial to initially}
Suppose $(i, c) \notin D$.
If $L_{i,c}$ only moves cells in the $(i+1,c)$-initial
segment of $D$, then it acts initially on $D$.
\end{lemma}
\begin{proof}
Suppose to the contrary that  $D$ is not fixed by $L_{i+1,c}$. 
Let $(r, c)$ be the first cell moved by $L_{i+1, c}$.
Clearly, $(r,c)$ is not in the 
$(i+1,c)$-initial segment of $D$ and
it will also be moved by $L_{i,c}$.
\end{proof}

We introduce more definitions that captures the structure of columns for intermediate diagrams during our algorithm. 

\begin{defn}
\label{D: paired}
We say a diagram $D$ is 
\definition{$(i,c)$-paired}
if the following are satisfied:
\begin{itemize}
\item Take any cell $(R,c) \in D$ with $i \leq R$ and $(R, c+1) \notin D$. 
There exists $(r, c+1) \in D$ 
with $i \leq r < R$ and $(r, c) \notin D$.
Moreover, $(r', c), (r', c+1) \in D$
for any $r < r' < R$.
\item Take any cell $(r, c+1) \in D$ 
with $i \leq r$ and $(r, c) \notin D$. 
There exists $(R,c) \in D$ with $r < R$ and $(R, c+1) \notin D$.
Moreover, $(r', c), (r', c+1) \in D$
for any $r < r' < R$.
\end{itemize}

\end{defn}

\begin{rem}
\label{R: paired and initially}
Notice that if $D$ is $(i,c)$-paired, 
then $L_{i,c}$ fixes $D$.
\end{rem}

\begin{exa}
Consider the following diagram $D$.
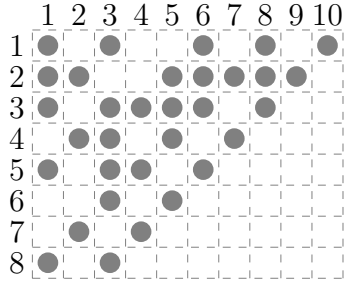
    $$
    \begin{tikzpicture}[x=1em,y=1em,thick,color = blue]
        \draw[step=1,gray,ultra thin,dashed] (0,0) grid (10,8);
        \node[color=black] at (0.5,8.5) {$1$};
        \node[color=black] at (1.5,8.5) {$2$};
        \node[color=black] at (2.5,8.5) {$3$};
        \node[color=black] at (3.5,8.5) {$4$};
        \node[color=black] at (4.5,8.5) {$5$};
        \node[color=black] at (5.5,8.5) {$6$};
        \node[color=black] at (6.5,8.5) {$7$};
        \node[color=black] at (7.5,8.5) {$8$};
        \node[color=black] at (8.5,8.5) {$9$};
        \node[color=black] at (9.5,8.5) {$10$};
        \node[color=black] at (-.5,7.5) {$1$};
        \node[color=black] at (-.5,6.5) {$2$};
        \node[color=black] at (-.5,5.5) {$3$};
        \node[color=black] at (-.5,4.5) {$4$};
        \node[color=black] at (-.5,3.5) {$5$};
        \node[color=black] at (-.5,2.5) {$6$};
        \node[color=black] at (-.5,1.5) {$7$};
        \node[color=black] at (-.5,0.5) {$8$};
        \filldraw [gray] (0.5,7.5) circle (3.5pt);
        \filldraw [gray] (0.5,6.5) circle (3.5pt);
        \filldraw [gray] (0.5,5.5) circle (3.5pt);
        \filldraw [gray] (0.5,3.5) circle (3.5pt);
        \filldraw [gray] (0.5,0.5) circle (3.5pt);
        \filldraw [gray] (1.5,6.5) circle (3.5pt);
        \filldraw [gray] (1.5,4.5) circle (3.5pt);
        \filldraw [gray] (1.5,1.5) circle (3.5pt);
        \filldraw [gray] (2.5,7.5) circle (3.5pt);
        \filldraw [gray] (2.5,5.5) circle (3.5pt);
        \filldraw [gray] (2.5,4.5) circle (3.5pt);
        \filldraw [gray] (2.5,3.5) circle (3.5pt);
        \filldraw [gray] (2.5,2.5) circle (3.5pt);
        \filldraw [gray] (2.5,0.5) circle (3.5pt);
        \filldraw [gray] (3.5,5.5) circle (3.5pt);
        \filldraw [gray] (3.5,3.5) circle (3.5pt);
        \filldraw [gray] (3.5,1.5) circle (3.5pt);
        \filldraw [gray] (4.5,6.5) circle (3.5pt);
        \filldraw [gray] (4.5,5.5) circle (3.5pt);
        \filldraw [gray] (4.5,4.5) circle (3.5pt);
        \filldraw [gray] (4.5,2.5) circle (3.5pt);
        \filldraw [gray] (5.5,7.5) circle (3.5pt);
        \filldraw [gray] (5.5,6.5) circle (3.5pt);
        \filldraw [gray] (5.5,5.5) circle (3.5pt);
        \filldraw [gray] (5.5,3.5) circle (3.5pt);
        \filldraw [gray] (6.5,6.5) circle (3.5pt);
        \filldraw [gray] (6.5,4.5) circle (3.5pt);
        \filldraw [gray] (7.5,7.5) circle (3.5pt);
        \filldraw [gray] (7.5,6.5) circle (3.5pt);
        \filldraw [gray] (7.5,5.5) circle (3.5pt);
        \filldraw [gray] (8.5,6.5) circle (3.5pt);
        \filldraw [gray] (9.5,7.5) circle (3.5pt);
    \end{tikzpicture}
    $$
Then $D$ has the following properties: $(1,5)$-paired, $(1,9)$-paired, $(4,1)$-paired, $(6,1)$-paired.
\end{exa}

We have the following lemma regarding this new notion. 

\begin{lemma}
    \label{L: paired induction}
    Let diagram $D$ be $(3,c)$-paired
    and $(2, c+1) \notin D$.
    We consider the action of $L_{1,c+1}L_{2,c}L_{3,c-1}$ on $D$. 
    Assume $L_{3, c-1}$ and $L_{2,c}$ act initially. 
    Let $(r_1, c), \cdots, (r_m,c)$
    be the cells moved by $L_{2,c}$
    with $r_1 < \cdots < r_m$
    and let $r_0 = 2$.
    We further assume $L_{1,c+1}$
    moves $(r_1', c+1), \cdots, (r_m',c+1)$
    with $r_{i-1} \leq r_i' < r_i$.
    Then $D' = L_{1,c+1}L_{2,c}L_{3,c-1}(D)$ is $(2,c)$-paired. 
\end{lemma}

\begin{exa}
    ~\label{E: paired}
    Consider the action of $L_{1,c+1}L_{2,c}L_{3,c-1}$ on $D$ whose column $c$ and $c+1$
    are depicted in the left-most figure. 
    We see $D$ is $(3,c)$-paired. 
    The action of $L_{2,c}$
    and $L_{1, c+1}$ satisfy
    the condition in Lemma~\ref{L: paired induction}: 
    For instance, $L_{2,c}$
    moves $(5,c)$ to $(2, c+1)$
    and there is a unique cell
    $(r,c+1)$ moved by $L_{1,c+1}$ with $2 \leq r < 5$, namely $(3, c+1)$.
    Then by the Lemma, 
    we know $L_{1,c+1}L_{2,c}L_{3,c-1}(D)$, whose column $c$ and $c+1$ are depicted in the right-most figure,
    is $(2,c)$-paired.
    
    $$
    
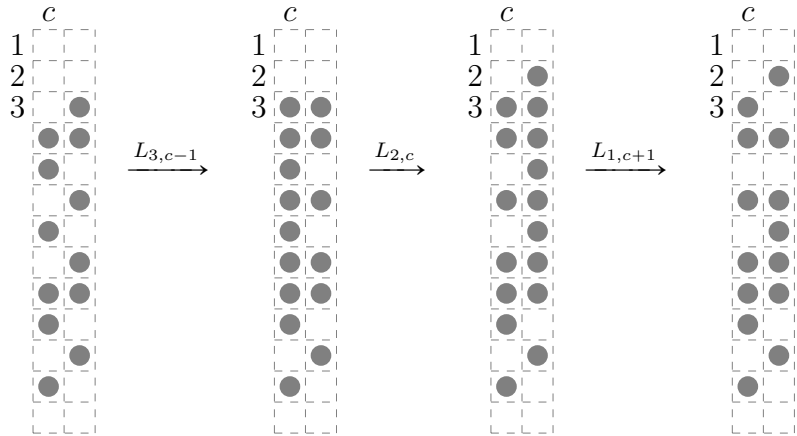
\begin{tikzpicture}[x=1em,y=1em,thick,color = blue]
        \draw[step=1,gray,ultra thin,dashed] (0,0) grid (2,13);
        \node[color=black] at (0.5,13.5) {$c$};
        \node[color=black] at (-0.5,12.5) {$1$};
        \node[color=black] at (-0.5,11.5) {$2$};
        \node[color=black] at (-0.5,10.5) {$3$};
        \filldraw [gray] (1.5,10.5) circle (3.5pt);
        \filldraw [gray] (0.5,9.5) circle (3.5pt);
        \filldraw [gray] (1.5,9.5) circle (3.5pt);
        \filldraw [gray] (0.5,8.5) circle (3.5pt);
        \filldraw [gray] (1.5,7.5) circle (3.5pt);
        \filldraw [gray] (0.5,6.5) circle (3.5pt);
        \filldraw [gray] (1.5,5.5) circle (3.5pt);
        \filldraw [gray] (0.5,4.5) circle (3.5pt);
        \filldraw [gray] (1.5,4.5) circle (3.5pt);
        \filldraw [gray] (0.5,3.5) circle (3.5pt);
        \filldraw [gray] (1.5,2.5) circle (3.5pt);
        \filldraw [gray] (0.5,1.5) circle (3.5pt);
    \end{tikzpicture}
    \quad \raisebox{3.4cm}{$\xrightarrow{L_{3,c-1}}$}\quad
    \begin{tikzpicture}[x=1em,y=1em,thick,color = blue]
        \draw[step=1,gray,ultra thin,dashed] (0,0) grid (2,13);
        \node[color=black] at (0.5,13.5) {$c$};
        \node[color=black] at (-0.5,12.5) {$1$};
        \node[color=black] at (-0.5,11.5) {$2$};
        \node[color=black] at (-0.5,10.5) {$3$};
        \filldraw [gray] (0.5,10.5) circle (3.5pt);
        \filldraw [gray] (1.5,10.5) circle (3.5pt);
        \filldraw [gray] (0.5,9.5) circle (3.5pt);
        \filldraw [gray] (1.5,9.5) circle (3.5pt);
        \filldraw [gray] (0.5,8.5) circle (3.5pt);
        \filldraw [gray] (0.5,7.5) circle (3.5pt);
        \filldraw [gray] (1.5,7.5) circle (3.5pt);
        \filldraw [gray] (0.5,6.5) circle (3.5pt);
        \filldraw [gray] (0.5,5.5) circle (3.5pt);
        \filldraw [gray] (1.5,5.5) circle (3.5pt);
        \filldraw [gray] (0.5,4.5) circle (3.5pt);
        \filldraw [gray] (1.5,4.5) circle (3.5pt);
        \filldraw [gray] (0.5,3.5) circle (3.5pt);
        \filldraw [gray] (1.5,2.5) circle (3.5pt);
        \filldraw [gray] (0.5,1.5) circle (3.5pt);
    \end{tikzpicture}
    \quad \raisebox{3.4cm}{$\xrightarrow{L_{2,c}}$}\quad
    \begin{tikzpicture}[x=1em,y=1em,thick,color = blue]
        \draw[step=1,gray,ultra thin,dashed] (0,0) grid (2,13);
        \node[color=black] at (0.5,13.5) {$c$};
        \node[color=black] at (-0.5,12.5) {$1$};
        \node[color=black] at (-0.5,11.5) {$2$};
        \node[color=black] at (-0.5,10.5) {$3$};
        \filldraw [gray] (1.5,11.5) circle (3.5pt);
        \filldraw [gray] (0.5,10.5) circle (3.5pt);
        \filldraw [gray] (1.5,10.5) circle (3.5pt);
        \filldraw [gray] (0.5,9.5) circle (3.5pt);
        \filldraw [gray] (1.5,9.5) circle (3.5pt);
        \filldraw [gray] (1.5,8.5) circle (3.5pt);
        \filldraw [gray] (0.5,7.5) circle (3.5pt);
        \filldraw [gray] (1.5,7.5) circle (3.5pt);
        \filldraw [gray] (1.5,6.5) circle (3.5pt);
        \filldraw [gray] (0.5,5.5) circle (3.5pt);
        \filldraw [gray] (1.5,5.5) circle (3.5pt);
        \filldraw [gray] (0.5,4.5) circle (3.5pt);
        \filldraw [gray] (1.5,4.5) circle (3.5pt);
        \filldraw [gray] (0.5,3.5) circle (3.5pt);
        \filldraw [gray] (1.5,2.5) circle (3.5pt);
        \filldraw [gray] (0.5,1.5) circle (3.5pt);
    \end{tikzpicture}
    \quad \raisebox{3.4cm}{$\xrightarrow{L_{1,c+1}}$}\quad
    \begin{tikzpicture}[x=1em,y=1em,thick,color = blue]
        \draw[step=1,gray,ultra thin,dashed] (0,0) grid (2,13);
        \node[color=black] at (0.5,13.5) {$c$};
        \node[color=black] at (-0.5,12.5) {$1$};
        \node[color=black] at (-0.5,11.5) {$2$};
        \node[color=black] at (-0.5,10.5) {$3$};
        \filldraw [gray] (1.5,11.5) circle (3.5pt);
        \filldraw [gray] (0.5,10.5) circle (3.5pt);
        \filldraw [gray] (0.5,9.5) circle (3.5pt);
        \filldraw [gray] (1.5,9.5) circle (3.5pt);
        \filldraw [gray] (0.5,7.5) circle (3.5pt);
        \filldraw [gray] (1.5,7.5) circle (3.5pt);
        \filldraw [gray] (1.5,6.5) circle (3.5pt);
        \filldraw [gray] (0.5,5.5) circle (3.5pt);
        \filldraw [gray] (1.5,5.5) circle (3.5pt);
        \filldraw [gray] (0.5,4.5) circle (3.5pt);
        \filldraw [gray] (1.5,4.5) circle (3.5pt);
        \filldraw [gray] (0.5,3.5) circle (3.5pt);
        \filldraw [gray] (1.5,2.5) circle (3.5pt);
        \filldraw [gray] (0.5,1.5) circle (3.5pt);
    \end{tikzpicture}
    $$
\end{exa}
\begin{proof}
Say $(t,c)$ is the bottom-most
cell in the $(2,c)$-initial segment of $L_{3, c-1}(D)$.
Since $L_{3,c-1}$ acts initially on $D$,
it will only move cells to the $(2,c)$-initial segment by Lemma~\ref{L: intoinitial}.
Since $L_{2, c}$ acts initially on $D$,
it will only move cells in the $(2,c)$-initial segment
by Lemma~\ref{L: movablecell}.
Then by our assumption in the lemma,
$L_{1,c+1}$ also moves cells above row $t$.
Thus, $D$ and $D'$ agreed under row $t$
in column $c$ and $c+1$.
Now we check $D'$ is $(2,c)$-paired.

Take $(R,c)$ in $D'$ such that $R \geq 2$
and $(R, c+1) \notin D'$.
We find the $r$ satisfying the condition 
in the definition of $(2,c)$-paired
by considering two cases.
\begin{itemize}
\item If $R > t$,
then $(R,c) \in D$ and $(R, c+1) \notin D$.
Since $D$ is $(3,c)$-paired, 
we can find $(r,c+1) \in D$
such that $2 \leq r < R$, $(r, c) \notin D$
and $(r', c), (r', c+1) \in D$ for $r < r' < R$.
It remains to show $r > t$.
If not, $(r, c)$ is in the $(2,c)$-initial segment
of $L_{3, c-1}(D)$, then so is $(R, c)$,
contradicting to $R> t$.
\item If $R \leq t$, then $(R, c) \in L_{3, c-1}(D)$.
If $(R, c+1) \notin L_{3, c-1}(D)$,
by Lemma~\ref{L: movablecell},
$L_{2,c}$ moves $(R,c)$.
Since $(R,c)$ is in $D'$,
we know it is the last cell moved by $L_{2,c}$,
so $R = r_m$.
By Lemma~\ref{L: chainmove}, $L_{2,c}$
moves $(r_m, c)$ to $(r_{m-1}, c+1)$.
We have $(r_{m-1}, c) \notin D'$.
By our assumption on $L_{1, c-1}$,
it does not make a regular ladder move on cells
between row $r_{m-1}$ and row $r_m$.
Thus, we may pick $r = r_{m-1}$.

Now assume $(R, c+1) \in L_{3, c-1}(D)$.
Then, $L_{1, c+1}$ moves $(R, c+1)$,
so $R = r_i'$ for some $i \in [m-1]$.
We know $L_{2,c}$ moves $(r_i, c)$
to $(r_{i-1}, c+1)$.
By our assumption on $L_{1,c+1}$,
$r_{i-1} < r_{i}'$ and 
$L_{1,c+1}$ does not make a move between
row $r_{i-1}$ and $r_i'$.
Thus, we may pick $r = r_{i-1}$.
\end{itemize}

Take $(r,c+1)$ in $D'$ such that $r \geq 2$
and $(r, c) \notin D'$.
We find the $R$ satisfying the condition 
in the definition of $(2,c)$-paired
by considering two cases.
\begin{itemize}
\item If $r > t$, then $(r,c+1) \in D$
and $(r,c) \notin D$.
Moreover, 
since $(2, c+1) \notin D$,
we know $r \geq 3$.
By $D$ is $(3,c)$-paired, 
we can find $R > r > t$
such that $(R,c) \in D$, $(R,c+1) \notin D$
and $(r', c), (r', c+1) \in D$ for $r < r' < R$.

\item If $r \leq t$, then $(r, c) \in L_{3, c-1}(D)$.
We know $L_{2,c}$ performs a regular 
ladder move on $(r,c)$, so $r = r_i$
for some $i \in [m-1]$.
We know $r_i < r_{i+1}' < r_{i+1}$
and $(r', c), (r', c+1) \in L_{2,c}L_{3,c-1}(D)$
for $r_i < r' < r_{i+1}$.
If $i + 1 < m$,
then $L_{1,c+1}$ makes a regular ladder move
on $(r_{i+1}', c+1)$.
We have $(r_{i+1}', c) \in D'$ and $(r_{i+1}, c+1) \notin D'$.
We may pick $R = r'$.
If $i+1 = m$,
then $L_{1,c+1}$ makes a K-ladder move
on $(r_{i+1}', c+1)$.
We may pick $R = r_{i+1}$.
\qedhere
\end{itemize}
\end{proof}

The last piece of our preparation work is the following observation.
\begin{rem}
\label{R: Swap}
Notice that $L_{i,c}$ and $L_{i', c'}$ commute
if $|c - c'| > 1$.
Therefore, we know applying
$$
L_{1,1} L_{1,2} \cdots L_{1,n-2} \quad L_{2,1} L_{2,2} \cdots L_{2,n-3}
$$
is the same as applying
\begin{align*}
L_{1,1} \:\: L_{1, 2} L_{2,1}  \:\: L_{1,3} L_{2,2} \:\: \cdots \:\: L_{1,n-4} L_{2,n-3} \:\: L_{1,n-2} L_{2,n-3}.
\end{align*}
Moreover, each $L_{i,c}$
behaves the same in both expressions. 
\end{rem}

Now we embark on proving Proposition~\ref{P: Movecode}
and Corollary~\ref{C: Row +1}.
We start by introducing two claims
which will imply Proposition~\ref{P: Movecode}
and Corollary~\ref{C: Row +1} respectively.
For a diagram $D$, let $D^{\downarrow k}$
be the diagram obtained by shifting 
all cells of $D$ downward by $k$.
We claim:
\begin{itemize}
\item Claim 1: Take $N \in \mathbb{Z}_{>0}$ and $w \in S_N$.
Consider
\begin{align}
\label{EQ: Claim 1}
(L_{1, 2} L_{2, 1}) \cdots (L_{1, N-2} L_{2, N-3})(L_{1, N} L_{2, N-1}) (\widehat{P}(w)^{\downarrow 2}).
\end{align}
Take any $c \in [N-1]$.
Then $L_{2, c}$ and $L_{1, c+1}$
moves the same number of cells. 
More specifically, 
suppose $L_{2, c}$ moves a cell 
$(r,c)$ to $(\hat{r},c+1)$.
Then there exists a unique $r'$
such that $\hat{r} \leq r' < r$
and $(r', c+1)$ is moved by  $L_{1, c+1}$.
In addition, 
after the action of $L_{1, c+1}$,
the diagram is $(2,c)$-paired.
\item Claim 2: Take $N \in \mathbb{Z}_{>0}$ and $w \in S_N$.
Consider
$$
L_{1, 1} \cdots L_{1, N-1} (\widehat{P}(w)^{\downarrow 1}).
$$
Each $L_{1, c}$ acts initially. 
\end{itemize}

We will inductively show both claims hold for all $N$.
The induction is based on Lemma~\ref{L: 1 to 2} 
and Lemma~\ref{L: 2 to 1}.

\begin{lemma}
\label{L: 1 to 2}
Suppose Claim 1 and Claim 2 hold for $N \leq n$,
then Claim 2 holds for $N = n + 1$.
\end{lemma}
\begin{proof}
Suppose $w = (b,u) \in S_{n+1}$. Let $D$ be the diagram obtained by putting $b$ left-justified cells in the second row
of $\widehat{P}(u)^{\downarrow 2}$. Then 
$\widehat{P}(w)^{\downarrow 1} = L_{2,1}L_{2,2} \cdots L_{2,n-1}(D)$
and each $L_{2,c}$ acts initially by Claim $2$ for $u$.
By Remark~\ref{R: Swap}, 
we may write $L_{1, 1} \cdots L_{1, N-1} (\widehat{P}(w)^{\downarrow 1})$
as 
\begin{align}
\label{EQ: Claim 2 proof}
L_{1, 1} \cdots L_{1, N-1} \: L_{2,1} \cdots L_{2,n-1}(D) = (L_{1, 2} L_{2, 1}) \cdots (L_{1, N-2} L_{2, N-3})(L_{1, N} L_{2, N-1}) (D).
\end{align}

Clearly, for $c \leq b$, $L_{1,c}$ acts initially on 
$\widehat{P}(w)^{\downarrow 1}$.
Now take $c> b$.
We know the $L_{1,c}$ behaves the same 
in both sides of~(\ref{EQ: Claim 2 proof}). 
By Lemma~\ref{L: initial to initially},
it is enough to show each $L_{1,c}$ on the right hand side moves cells in the $(2,c)$-initial segment.
Since $L_{2,c-1}$ acts initially, by Lemma~\ref{L: intoinitial}, $L_{2,c-1}$ move cells into the $(2,c)$-initial segment. 
Then by claim $1$ of $u$, $L_{1,c}$ 
moves cells in the $(2,c)$-initial segment. 
\end{proof}

\begin{lemma}
\label{L: 2 to 1}
Suppose Claim 1 holds for $N \leq n$
and Claim 2 holds for $N \leq n+1$,
then Claim 1 holds for $N = n + 1$.
\end{lemma}
\begin{proof}
Since Claim 2 holds for $N \leq n+1$, each $L_{1,c}$ and $L_{2,c}$ 
in~(\ref{EQ: Claim 1}) acts initially by Remark~\ref{R: Swap}.
We prove Claim 1 by induction on $c = n, \cdots, 2,1$.
The base case with $c = n$ is trivial. 

Suppose $c \in [n-1]$. 
Let $D'$ be the diagram right before applying $L_{2,c}$ in~(\ref{EQ: Claim 1}).
By our inductive hypothesis for $c + 1$,
$D'$ is $(2,c+1)$-paired.
Now apply $L_{2,c}$ to $D'$.
Let $(r_1, c), \cdots, (r_k, c)$ be the cells
moved by $L_{2,c}$.
Let $r_0 = 2$.
For $j \in [k]$, by Lemma~\ref{L: chainmove},
$(r_j,c)$ is moved to $(r_{j-1},c+1)$.
By Lemma~\ref{L: intoinitial}, 
$(r_{j-1},c+1)$ is in the $(2, c+1)$-initial segment of 
$L_{2,c}(D)$.
We consider two cases.
\begin{itemize}
\item If $(r_{j-1},c+2) \notin D'$, 
then $(r_{j-1},c+1)$ will be moved by $L_{1,c+1}$ by Lemma~\ref{L: movablecell}. 
For $r_{j-1} < r' < r$,
by $D'$ is $(2, c+1)$-paired,
we know $(r', c+1), (r', c+2) \in D'$.
By Lemma~\ref{L: movablecell},
$L_{1,c+1}$ will not move $(r', c+1)$.
\item Now assume $(r_{j-1},c+2) \in D'$.
Since $D'$ is $(2,c+1)$-paired
and $(r_{j-1},c+) \notin D'$,
we can find $R > r_{j-1}$ 
such that $(R,c+1) \in D'$, $(R,c+2) \notin D$
and $(r', c+1), (r', c+2) \in D'$
for any $r_{j-1} < r' < R$. 
We know $(r_j, c+1) \notin D'$,
so $R < r_j$.
For $R < r' < r_j$, 
since $(r', c+1) \in D'$ 
and $D'$ is $(2, c+1)$-paired,
we must have $(r', c+2) \in D'$.
By~\ref{L: movablecell}, $(R ,c+1)$ is the unique cell moved during $L_{1,c+1}$ between row $r_{j-1}$ and row $r_j$.
\end{itemize}

Now we show $L_{1, c+1}$ and $L_{2,c}$ 
move the same number of cells, 
we already know $L_{1,c+1}$ makes exactly
one move between
row $r_{j-1}$ and row $r_j$ inclusively
for $j \in [k]$. 
We just need to show $L_{1,c+1}$
does not move any $(r,c+1)$ for any $r > r_k$.
Notice that $(r_k, c+1) \notin L_{2,c}(D')$,
so $(r, c+1)$ is not in the $(2, c+1)$-initial
segment of $L_{2,c}(D')$.
By Lemma~\ref{L: movablecell},
$(r, c+1)$ will not be moved.

It remains to check $L_{1,c+1} L_{2,c}(D')$ is
$(2,c)$-paired.
Write $w$ as $(b,u)$.
Let $D$ be the diagram obtained by putting $b$ left-justified cells in row $3$
of $\widehat{P}(u)^{\downarrow 3}$. Then 
\begin{align*}
    \widehat{P}(w)^{\downarrow 2} = L_{3,1}L_{3,2} \cdots L_{3,n-1}(D)
\end{align*}
By Remark~\ref{R: Swap}, 
\begin{align*}
    &(L_{1, 2} L_{2, 1}) \cdots (L_{1, n+1} L_{2, n}) (\widehat{P}(w)^{\downarrow 2})\\
    = & (L_{1, 2} L_{2, 1}) \cdots (L_{1, n+1} L_{2, n}) (L_{3,1}L_{3,2} \cdots L_{3,n-1})(D)\\
    = & (L_{1,2} L_{2,1})(L_{1, 3} L_{2, 2} L_{3,1}) \cdots (L_{1, n+1} L_{2, n} L_{3,n-1})(D)
\end{align*}

If $c > b$,
then $(3,c) \notin D$.
By claim $1$ of $u$, after $L_{2,c+1}$ the diagram is $(3,c)$-paired. Therefore, by Lemma~\ref{L: paired induction}, after $L_{1,c+1}$ the diagram is $(2,c)$-paired. 

Now consider $c \leq b$, so $(3,c) \in D$. 
We consider three cases:
\begin{itemize}
\item Case 1: $(3,c)$ is moved by $L_{2,c}$ and not the last cell moved by $L_{2,c}$.
Then $L_{2,c}$ performs a regular ladder move
on $(3,c)$ moving it to $(2, c+1)$.
Later, $L_{1,c+1}$ will move $(2,c+1)$.
Since $L_{1,c+1}$ and $L_{2,c}$ moves the same number of cells,
we know $L_{1,c+1}$ makes a regular ladder move on $(2, c+1)$.
By Remark~\ref{R: L recursion},
the action of $L_{1, c+1} L_{2,c}$
is the same as first moving $(3,c)$ to $(1, c+2)$,
and then perform $L_{2, c+1}L_{1, c+2}$.
By Claim $1$ of $u$,
the diagram after applying $L_{1,c+1}$ is $(3,c)$-paired.
Since $(2,c), (2,c+1)$ are not in the diagram, 
it is $(2,c)$-paired. 

\item 
Case 2: $(3,c)$ is the last cell moved by $L_{2,c}$.
Then $L_{2,c}$ performs a K-ladder move
on $(3,c)$ moving it to $(2, c+1)$.
Later, $L_{1,c+1}$ will move $(2,c+1)$.
Since $L_{1,c+1}$ and $L_{2,c}$ moves the same number of cells,
we know $L_{1,c+1}$ makes K-ladder move on $(2, c+1)$.
By Remark~\ref{R: L recursion},
the action of  $L_{1, c+1} L_{2,c}$
can be described as follows: 
Remove $(3,c)$, perform $L_{2, c+2}L_{3,c}$,
and then add cells $(3,c)$, $(2,c+1)$ and $(1, c+2)$.
By Claim $1$ of $u$,
before adding those three cells,
the diagram is $(3,c)$-paired.
Thus, after adding these three cells,
the diagram is $(2,c)$-paired.

\item If $(3,c)$ is not moved by $L_{2,c}$, then $(3,c+1) \in D$.
By Remark~\ref{R: L recursion}, applying $L_{1,c+1}L_{2,c}$ 
is the same as applying $L_{2, c+2}L_{3,c}$ while
ignoring row $3$.
By Claim $1$ of $u$,
after the action of $L_{1, c+1}$, the diagram is $(2,c)$-paired. \qedhere
\end{itemize}
\end{proof}

\begin{lemma}
\label{L: 1 and 2}
Claim $1$ and $2$ hold for all $N \in \mathbb{Z}_{>0}$.
\end{lemma}
\begin{proof}
The claims are obvious when $N = 1$.
Then we prove by induction on $N$.
The inductive step is given by Lemma~\ref{L: 1 to 2}
and Lemma~\ref{L: 2 to 1}.
\end{proof}

\begin{cor}
\label{C: Initially}
In~(\ref{EQ: algorithm}),
each $L_{i,c}$ acts initially. 
\end{cor}
\begin{proof}
Suppose $w\in S_n$ and we prove the corollary by
induction on $n$.
Suppose $w = (b,u)$.
Since the corollary holds for $u$,
we know $L_{i,c}$ in~(\ref{EQ: algorithm})
acts initially when $i > 2$.
Finally, each $L_{1,c}$
acts initially by Claim 2. 
\end{proof}

Now we may prove the main results of this subsection
using the two claims.
\begin{proof}[Proof of Proposition~\ref{P: Movecode}]
We induct on $n$. The base cases $n = 2$ is trivial. 
Now suppose $n > 2$ and take $v=(a,w) \in S_n$. 
Let $D$ be the diagram obtained by putting $a$ left-justified cells 
in row $1$ of $\widehat{P}(w)^{\downarrow 1}$.
The last iteration to compute $\widehat{P}(v)$ is to apply 
$L_{1,1} \cdots L_{1,n-2}L_{1,n-2}$ on $D$.
For $c \in [a]$,
since $L_{1,c}$ acts initially and $(1,c) \in D$,
$L_{1,c}$ does not move any cells.

Now assume $c > a$.
We want to show $L_{1,c}$ moves exactly $\movecode(w)_c$ cells. 
Let $w = (b,u)$ and let $D'$ be the diagram obtained by putting $b$ left-justified cells in the row $2$ of $\widehat{P}(u)^{\downarrow 2}$. Then, 
\begin{align*}
    &L_{1,1} \cdots L_{1,n-2}L_{1,n-1}(D) \\
    = &(L_{1,1} \cdots L_{1,n-2}L_{1,n-1})(L_{2,1} \cdots L_{2,n-3}L_{2,n-2})(D')\\
    = &(L_{1,1})(L_{1,2}L_{2,1})\cdots(L_{1,n-1}L_{2,n-2})(D')
\end{align*}

For $c>b$, by our induction hypothesis, applying $L_{2,c}$ moves exactly $\movecode(u)_c$ cells. Then by Claim $1$, applying $L_{1,c+1}$ to $D$ also moves exactly $\movecode(u)_c$ cells. Therefore the number of cells moved by $L_{1,c+1}$ is $\movecode(u)_c = \movecode(w)_{c+1}$.
Now clearly each $L_{2,c}$ does not move any cells for $c \in [b]$.
We know $L_{1,b+1}$ also moves no cells since the $(2,b+1)$-initial segment is empty. Therefore $L_{1,b+1}$ moves $0 = \movecode(w)_{b+1} $ cells. 

Let $c_0$ be the largest in $[b]$ such that $\movecode(u)_{c_0} = 0$.
Say $c_0 = 0$ if no such $c_0$ exists.
For $c \in [b]$,
by Lemma~\ref{L: constructmovecode},
we have
$$
\movecode(w)_c = \begin{cases}
\movecode(u)_c + 1 & \text{if $c \geq c_0$.}\\
\movecode(u)_c & \text{otherwise}
\end{cases}$$

We first inductively show that for $c = b,\cdots, c_0 + 1$,
there is no cell at $(2, c+1)$ 
right before the action of $L_{1, c}$,
so $L_{1,c}$ moves $(2,c)$.
Moreover, 
$L_{1,c}$ moves $\movecode(w)_c > 2$ cells,
so the move on $(2,c)$ is a regular ladder move.
For $c = b$,
we know $(2,b+1)$ is always empty. 
For $c_0 < c < b$,
we know $L_{1,c+1}$ makes a regular ladder move on $(2,c+1)$,
so $(2, c+1)$ is empty right before the action of $L_{1,c}$.
Now for $c = b,\cdots, c_0 + 1$,
after $L_{1,c}$ moves $(2,c)$,
it behaves as if $L_{2,c}$ by Remark~\ref{R: L recursion}.
Thus, the total number of cells moved
is $\movecode(u)_c + 1 = \movecode(w)_c$.

Now consider $L_{1,c_0}$
when $c_0 > 0$.
Right before its action,
$(2,c_0+1)$ is empty.
Thus, $L_{1,c_0}$ will first move $(2,c_0)$ 
to $(1,c_0+1)$.
After that,
the number of cells it moves
is $\movecode(u)_{c_0}$, 
which is zero. 
Thus, the move on $(2,c_0)$ 
is a K-ladder move.
Also, $L_{1,c_0}$ moves $1 = \movecode(w)_{c_0}$ cell.

Finally, we prove by induction that
for $c = c_0 - 1, \cdots, 1$,
right before the action of $L_{1,c}$, 
the diagram contains $(2,c)$ and $(2,c+1)$.
For the base case, 
right before the action of $L_{1, c_0-1}$,
we know $(2,c_0)$ is in the diagram.
Now assume 
right before the action of $L_{1,c}$, 
the diagram contains $(2,c)$ and $(2,c+1)$ for some $c < c_0$.
Then $L_{1,c}$ will not move $(2,c)$.
After the action of $L_{1,c}$,
we know $(2,c)$ is still in the diagram. The inductive step is finished.
Now by Remark~\ref{R: L recursion},
the action of $L_{1,c}$ moves 
the same number of cells 
as $L_{2,c}$ on the diagram
without 
$(2,c)$ and $(2,c+1)$.
Thus, $L_{1,c}$ makes
$\movecode(u)_c = \movecode(w)_c$
moves. \qedhere
\end{proof}

\begin{proof}[Proof of Corollary~\ref{C: Row +1}]
Implied by Corollary~\ref{C: Initially} and Lemma~\ref{L: chainmove}.
\end{proof}

\section{Acknowledgements}
We thank Jianping Pan and Brendon Rhoades for
carefully reading an earlier version of this paper and giving many useful comments.
We are especially grateful to Anna Weigandt for suggesting this problem.
This project was performed as a Research Experience for Undergraduates at UC San Diego in summer 2023.

\bibliographystyle{alpha}
\bibliography{main.bbl}{}
\end{document}